\documentclass[preprint]{elsarticle}
\makeatletter{}\makeatletter{}\usepackage{floatrow}
\usepackage[margin=1.5in]{geometry}

\usepackage{color}
\usepackage[usenames,dvipsnames]{xcolor}
\usepackage[export]{adjustbox}

\usepackage[font=footnotesize]{subcaption}
\usepackage{morefloats}
\usepackage{booktabs}
\usepackage{multirow}
\usepackage{tablefootnote}
\usepackage[section]{placeins}

\usepackage{microtype}
\usepackage{siunitx}
\usepackage[super]{nth}

\usepackage[pdftex,colorlinks,allcolors=black]{hyperref}
\usepackage{xspace}
\usepackage{amsxtra}

\usepackage[capitalise,noabbrev]{cleveref}
\crefformat{appendix}{#2#1#3}
\Crefformat{equation}{Equation~(#2#1#3)}
\crefformat{equation}{(#2#1#3)}

\makeatletter{}

\renewcommand{\vec}[1]{\mathbf{#1}}

\makeatletter{}\newcommand{\packagefont}[1]{{#1}\xspace}
\newcommand{\fenics}{\packagefont{FEniCS}}
\newcommand{\otf}{\packagefont{OpenTidalFarm}}
\newcommand{\dolfin}{\packagefont{DOLFIN}}
\newcommand{\scipy}{\packagefont{SciPy}}
\newcommand{\da}{\packagefont{dolfin-adjoint}}

\makeatletter
\def\ps@pprintTitle{ \let\@oddhead\@empty
 \let\@evenhead\@empty
 \def\@oddfoot{} \let\@evenfoot\@oddfoot}
\makeatother

\begin{document}
\begin{frontmatter}

\title{Hybrid global-local optimisation algorithms for the layout design of
  tidal turbine arrays}
\author{G.L.~Barnett\corref{cor1}}
\ead{george.barnett10@imperial.ac.uk}
\author{S.W.~Funke}
\author{M.D.~Piggott}

\cortext[cor1]{Corresponding author}

\begin{keyword}
marine renewable energy \sep
tidal turbines \sep
gradient-based optimisation \sep
genetic algorithm based optimisation \sep
adjoint method
\end{keyword}

\makeatletter{}\begin{abstract}
Tidal stream power generation represents a promising
source of renewable energy. In order to extract an
economically useful amount of power, tens to hundreds of tidal
turbines need to be placed within an array. The
layout of these turbines can have a significant impact on
the power extracted and hence on the viability of the site.
\citet{funke2014} formulated the question
of the best turbine layout as an optimisation problem
constrained by the shallow water equations and
solved it using a local, gradient-based optimisation algorithm.
Given the local nature of this approach, the question
arises of how optimal the layouts actually are.
This becomes particularly important for scenarios
with complex bathymetry and layout constraints, both of
which typically introduce locally optimal layouts.
Optimisation algorithms
which find the global optima generally require orders
of magnitude more iterations than local optimisation
algorithms and are thus infeasible in combination with an expensive
flow model. This paper presents
an analytical wake model to act as an efficient proxy to the shallow water
model. Based upon this, a hybrid
global-local two-stage optimisation approach is
presented in which turbine layouts are first optimised with
the analytical wake model via a global optimisation
algorithm, and then further optimised with the shallow
water model via a local gradient-based optimisation algorithm.
This procedure is applied to
a number of idealised cases and a more realistic case
with complex bathymetry in the Inner Sound of the Pentland
Firth, Scotland. It is shown that in cases where bathymetry
is considered, the two-stage optimisation procedure is able
to improve the power extracted from the array by as much as
\SI{25}{\percent} compared to local optimisation for
idealised scenarios and by as much as \SI{12}{\percent} for
the more realistic Pentland Firth scenario whilst in many
cases reducing the overall computation time by approximately
\SIrange{30}{40}{\percent}.
\end{abstract}

\end{frontmatter}

\makeatletter{}\section{Introduction}\label{sec:introduction}

Tidal stream energy converters are one of the most promising technologies for
large-scale renewable power generation.  The use of marine turbines to
generate power has a key advantage over wind turbines: the stability and
predictability of the tidal currents which drive them. One of the greatest
challenges in the tidal renewable energy industry is in deciding upon the
precise location and arrangement of large arrays of tidal turbines, that is
finding the optimal location for the turbines in the array,
as it can have a significant effect on the total power and thus economic
performance and viability of a given site \citep{funke2014}.
Previous work on tidal turbine array optimisation can be broadly categorised
by the complexity of the flow model which controls the interaction between the
turbines and the driving currents. Simple flow models describe a reduction
in flow velocity due to the energy extracted by the turbines and are often
defined by analytical expressions or by expressions which are cheap to
compute. Optimal layouts are then either derived analytically or via an
extensive search of the parameter space using schemes such as genetic
algorithms (introduced in \cref{subsec:genetic}). These flow models typically
do not capture the complex nature of the flow interactions between turbines,
nor the potential for large-scale redirection of the resource due to blockage
effects, and feature more as estimates of the maximum energy that may be
extracted from a given site \citep{bryden2007,garrett2008}.

More complex models are usually formulated as numerical solutions to
partial differential equations (PDE). These models are
computationally expensive and an extensive search of the
parameter space is prohibitive.  Thus layout optimisation is often
performed manually, guided by intuition and experience
\citep{divett2013}. In realistic domains this task becomes difficult due to
the effects of bathymetry, complex flow patterns and large numbers of
turbines. \citet{funke2014} aim to alleviate this problem by posing
the layout problem as a PDE-constrained optimisation problem based on
the shallow water equations (SWE)
whilst leveraging the adjoint technique to efficiently solve the problem using
gradient-based optimisation algorithms. This approach is implemented
in the open-source framework \otf (\url{http://opentidalfarm.org}).

The adjoint technique is vital to making this approach viable as it allows to
compute the gradient in similar time to that of one model solve. Crucially,
for large arrays, this time remains almost independent of the number of
turbines to be optimised, which makes it feasible to optimise hundreds of
turbines. The layout optimisation consists of a number of steps. First, the
total power production of the array (the functional of interest) and its
gradient with respect to the position of each turbine in the array (the input
parameters) are evaluated. This information is then passed to a gradient-based
optimisation algorithm to relocate the turbines. This procedure is repeated
until an exit criteria is met.

\begin{figure*}[tb]
  \makebox[0pt]{    \includegraphics[center]{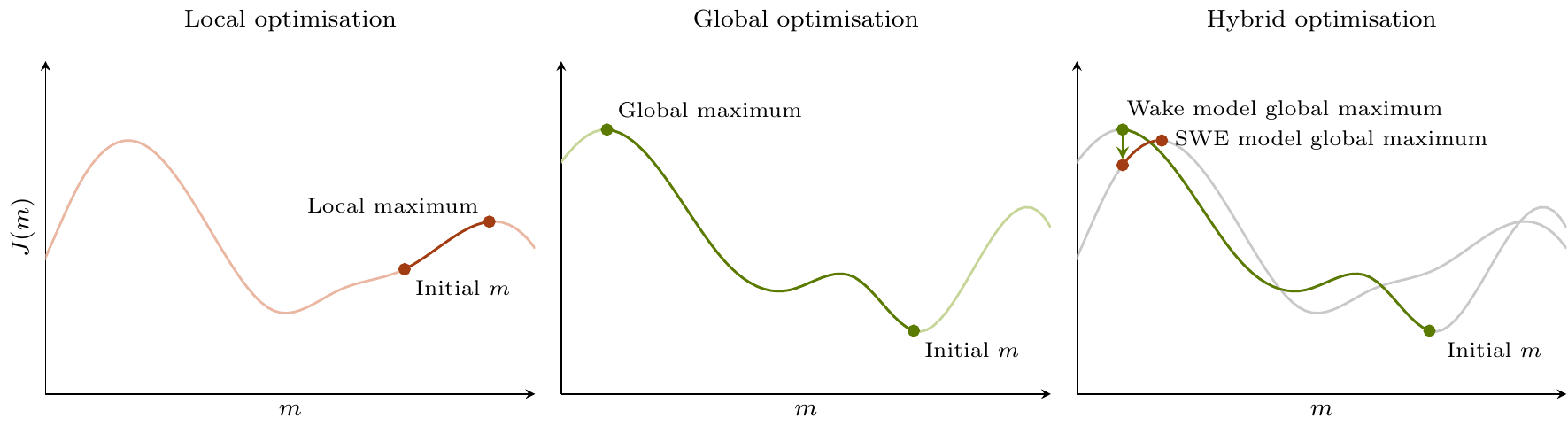}
  }
  \caption{Schematics showing local, global and hybrid optimisation of a
    functional of interest $J(m)$ with respect to its input parameter $m$.  The left
    schematic highlights the issue of multiple local maxima in the turbine
    optimisation problem whilst the central image highlights the utility of
    global optimisation. The proposed global-local hybrid approach is
    presented in the rightmost pane where the shallow water model is locally optimised
    from a good initial guess provided by a global optimisation of the wake
    model. The good initial guess provides a better chance of arriving at a
    better or even global maximum.
  }
  \label{fig:schematic}
\end{figure*}

The turbine layout problem is a nonlinear optimisation problem, and
can hence yield multiple local maxima, see \cref{fig:schematic}.
However, finding the global maximum with a global optimisation scheme
typically requires orders of magnitude more iterations than a local
optimisation method \citep{bilbao2009,chelouah2000}. Thus a direct
application of global methods on an expensive flow model is
infeasible, which is why \citet{funke2014} used local optimisation algorithms.

The work presented here expands the \otf framework to alleviate the problem of
combining a complex flow model with global optimisation algorithms by
providing an inexpensive, analytical model (herein referred to as the wake
model) to act as an efficient proxy to the complex flow model. A two-stage
approach to optimisation is then considered in the follow schemes:
\begin{enumerate}
  \item global optimisation of the wake model via a genetic algorithm followed
        by local gradient-based optimisation of the SWE model,
  \item global optimisation of the wake model using a basin-hopping
        gradient-based algorithm (introduced in \cref{subsec:gradient-based})
        followed by local gradient-based optimisation of the SWE model.
\end{enumerate}

The remaining sections are organised as follows.
\cref{sec:current-model} details the SWE model currently used in \otf,
analytical wake models used in the wind turbine industry and the wake model
developed for this work. An overview of the optimisation schemes used in this
research are discussed in \cref{sec:optimisation} with results from several
idealised scenarios presented and discussed in \cref{sec:results}. A more
realistic optimisation of 64 turbines in the Inner Sound of the Pentland
Firth, Scotland is presented in \cref{sec:orkney}. Finally, conclusions are
drawn in \cref{sec:conclusion}.

\makeatletter{}\section{Flow models}\label{sec:current-model}
\makeatletter{}\subsection{The shallow water equations in \otf}
\label{subsec:opentidalfarm}

The physical laws governing the tidal currents in \otf are modelled using the
two-dimensional nonlinear shallow water equations with individual
turbines being parameterised by regions of smoothly increasing friction
coefficients \citep{funke2014}. A similar approach is taken by
\citet{divett2013}. For simplicity, the steady equations are considered here:
\begin{equation}\label{eq:shallow-water}
  \begin{split}
    \vec{u}\cdot\nabla \vec{u} - \nu \nabla^2 \vec{u} + g \nabla \eta +
      \frac{c_b + c_t(\vec{m})}{H} \|\vec{u}\|\vec{u} = 0, \\
    \nabla \cdot \left(H\vec{u}\right) = 0,
  \end{split}
\end{equation}
where $\vec{u}$ is the depth-averaged horizontal velocity, $\nu$ is the
kinematic viscosity, $g$ is acceleration due to gravity, $\eta$ is
free-surface displacement, $H$ is the water depth at rest, whilst $c_b$ is the
drag coefficient for the natural bottom friction and
$c_t$ represents the additional friction induced by the turbine
parameterisation which depends on the turbine coordinates.  The coordinates
for $N$ turbines are stored in a vector $\vec{m}$ as:
\begin{equation}\label{eq:m}
  \vec{m} = \left(x_1, y_1, x_2, y_2, \ldots, x_{N}, y_{N} \right).
\end{equation}

\begin{figure}[tb]
  \centering
  \includegraphics[center]{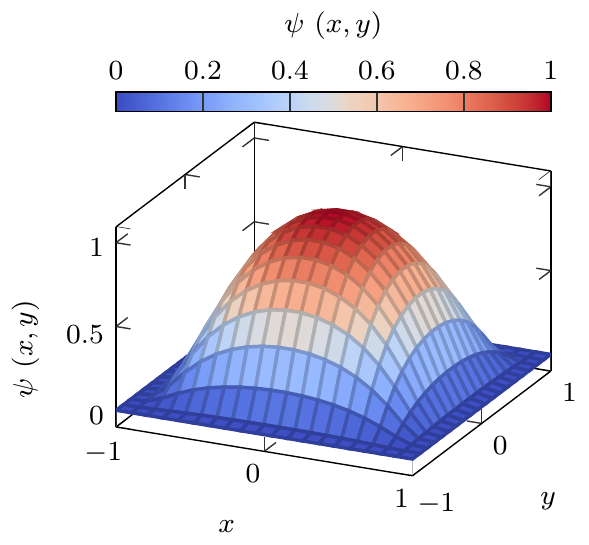}
  \caption{The turbines are parametrised by increasing friction in the turbine
    area in the shape of a bump function. The figure shows the bump function
    $\psi_{0,0,1}(x,y)$ centred on $(0,0)$ with a support radius of $1$ in
    both the $x$ and $y$ directions. This work relies on the fact that an
    individual turbine region is resolved by the mesh.}
  \label{fig:bump}
\end{figure}

The magnitude and extent of the added friction depends on the turbine of
interest. In order to be able to apply gradient-based optimisation, the
friction function $c_t(\vec{m})$ must be continuously differentiable with
respect to the turbine positions. Therefore, the parameterisation is chosen to
be a bump function $\psi(x,y)$, that is a smooth function with compact
support. This is derived in \citet{funke2014} and depicted in \cref{fig:bump}.

The friction function for a single turbine centred at $(x_i, y_i)$ with
support radius $r$ and a maximum friction value of $K_i$ is thus defined as:
\begin{equation}\label{eq:friction-function}
  C_i(\vec{m})(x,y) = K_i \psi_{x_i,y_i,r}(x, y).
\end{equation}
Finally, the dimensionless turbine friction function $c_t$ in
\cref{eq:shallow-water} is the sum of the bump functions for all $N$ turbines:
\begin{equation}\label{eq:friction-function-sum}
  c_t(\vec{m}) = \sum_{i=1}^{N} C_i(\vec{m}).
\end{equation}

The functional of interest $J$ describes the value to be
maximised and must be differentiable for gradient-based optimisation.
\citet{funke2014} chose the time-averaged power extracted from the tidal array
due to the increased bottom friction from the presence of turbines.
For the steady-state problem considered here the power extraction may be
expressed as:
\begin{equation}\label{eq:functional-of-interest}
  J(\vec{u}, \vec{m}) = \frac{1}{2}\int_\Omega \rho
    c_t(\vec{m})\|\vec{u}\|^3,
\end{equation}
where $\rho = \SI{1000}{\kg\per\m\cubed}$ is the fluid density  and
$\Omega$ is the domain of interest. For simplicity a cut-in speed and power
rating of the turbines is not considered.

\makeatletter{}\subsection{Analytical wake models}\label{subsec:wake-models}

From the functional of interest \cref{eq:functional-of-interest} it is clear
that power extracted from the tidal flow depends on the positions of the
turbines within the array, the flow field, and the applied friction due to the
presence of each turbine.

As currents interact with a turbine, energy is extracted, changing the
dynamics of the flow. Notably the flow is decelerated behind the turbine and
accelerated along its sides. The region around the turbine where flow
properties are changed significantly is the wake of the turbine. A simplified
model may base a computation of the extracted power of an array on the cube of
the flow speed \citep{vennell2012,funke2014} at each of the turbines. Thus the
power extracted from the array may be defined as:
\begin{equation}\label{eq:initial-power}
  J(\vec{u}, \vec{m}) = \alpha \sum_{i=1}^{N} \lVert \vec{u}_r(x_{i}, y_{i}) \rVert^3,
\end{equation}
where $\vec{u}_r$ is the `reduced' flow field due to the presence of the
wakes, and $\alpha$ is a constant incorporating fluid density, the turbine
specifications and can potentially also depend on the flow speed.

Simplified models do not include the effects of the turbines on the velocity
directly, and thus the dependence of the flow field on the turbine positions
must be carefully considered. The wake behind a turbine can be incorporated as
an analytical velocity deficit which decreases with increasing distance behind
the turbine.  Schemes similar to this have been used to optimise turbine
layouts within the wind turbine industry \citep{bilbao2009,changshui2011}. In
this case the reduced velocity field can be defined as:
\begin{equation}\label{eq:reduced-velocity}
  \vec{u}_r(x,y) = c(x,y;\vec{m}) \vec{u}_a(x,y),
\end{equation}
where $\vec{u}_r$ and $\vec{u}_a$ are the reduced and ambient flow fields
respectively, whilst $c$ is a combined velocity reduction factor field due to
the individual turbine reduction factors resulting from the wakes of the other
turbines in the array (discussed in \cref{subsec:wind-wake-models}).

It should be noted that analytical wake models of this type differ greatly
from a fully coupled model as described in \cref{subsec:opentidalfarm} in that
they only affect the local flow field around each turbine and do not impact
the flow field upstream of the turbine array, nor in the far-field.  In
particular, these simplifications preclude the conservation of energy, that is
the ambient flow is unaffected by energy extraction by turbines.

\makeatletter{}\subsection{Wind turbine wake models}\label{subsec:wind-wake-models}

A number of wake models exist for wind turbines which have been used for
optimising arrays of wind turbines. Of note are the Jensen and Larsen models
\citep{jensen1983,larsen1996}. Whilst these wake models are adequate for small
wind arrays \citep{barthelmie2007} they can under-estimate the wake losses in
large multi-row arrays \citep{mechali2006}.  \Citet{palm2010} consider such
models for the design of tidal turbine arrays; through comparison to CFD
simulations, wake models were concluded to hold promise for the appraisal of
tidal turbine arrays.

The Jensen model \citep{jensen1983} assumes a linearly expanding wake with a
velocity deficit that depends only on the distance behind the turbine.  The
model is defined by the diameter of its wake which depends on the turbine
thrust coefficient, a wake decay parameter and the diameter of the turbine.
An example wake from the Jensen model is depicted in
\cref{fig:wind-and-swe-wake}.

The Larsen model \cite{larsen1996} derives from Prandtl's turbulent boundary
layer equations. Both first- and second-order approximate solutions to the
boundary layer equations exist. It has been demonstrated that in practice
there is little difference between the solutions as discrepancies only lie in
the near wake field \citep{sorensen2008} where it would anyway be inefficient to
place another turbine. The diameter of the wake in the Larsen model derives
from the Prandtl mixing length and the rotor disc area. A flow field around a
turbine using the Larsen model is depicted in \cref{fig:wind-and-swe-wake}.

\begin{figure}[tb]
  \includegraphics[center]{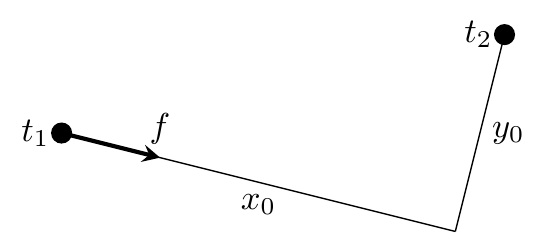}
  \caption{The relative $x$ and $y$ coordinates ($x_0, y_0$) of turbine $t_2$ in
           the coordinate system where $t_1$ is at the origin and the $x$-axis
           is aligned with $f$, the flow direction vector at $t_1$.}
  \label{fig:x0y0}
\end{figure}

However, both models presented fail to capture the flow fields
observed around tidal turbines. As both models are designed for wind turbines
they assume a free upper boundary and a fixed lower boundary --- the ground.
However, in tidal settings the turbine lies between the seabed (a fixed
boundary) and the free surface of the water. As there are constraints on the
depth in which a turbine may be placed it is likely that the existence of the
two bounding surfaces will act to restrict vertical wake expansion. This issue
is considered in the context of tidal turbines by \citet{thomson2011} where an
eddy-viscosity model \citep{ainslie1988} is solved locally to each turbine for
its downstream velocity deficit using a finite-difference method. The presence
of the bounding surfaces is accounted for by limiting vertical wake expansion
and extending lateral expansion through the use of an elliptical Gaussian
function. This approach is more complex and computationally expensive than the
completely analytical wake approaches considered here but may be considered in
future research.

A number of other issues exist with the wind wake models. Notably, the
parameters controlling the Larsen model are coupled such that it was not
possible to produce a wake with similar features to that of the wake behind a
marine turbine. Another numerical problem with the Jensen
model is that the derivative of the wake reduction factor in the direction
perpendicular to the ambient flow is zero. Thus a local gradient-based
optimisation algorithm will move turbines apart in the direction of the flow
before moving them in the direction perpendicular to the flow. Generally,
both the Larsen and Jensen models are problematic for gradient-based optimisation as the wake
reduction factor is discontinuous at the edge of the wake (excepting the
trailing edge). This has again been observed in numerical experiments where
turbines have been badly repositioned by the optimisation algorithm due to a
discontinuous gradient.

The core problem with both models is that they do not include zones where the
flow is accelerated locally. This feature of the wake is partially responsible
for the layouts featuring staggered `barrages' of turbines as demonstrated by
\citet{funke2014}, thus another wake model must be developed.

\makeatletter{}\subsection{Approximate shallow water wake model}\label{subsec:swe-wake-model}

\begin{figure*}[tb]
  \includegraphics[center]{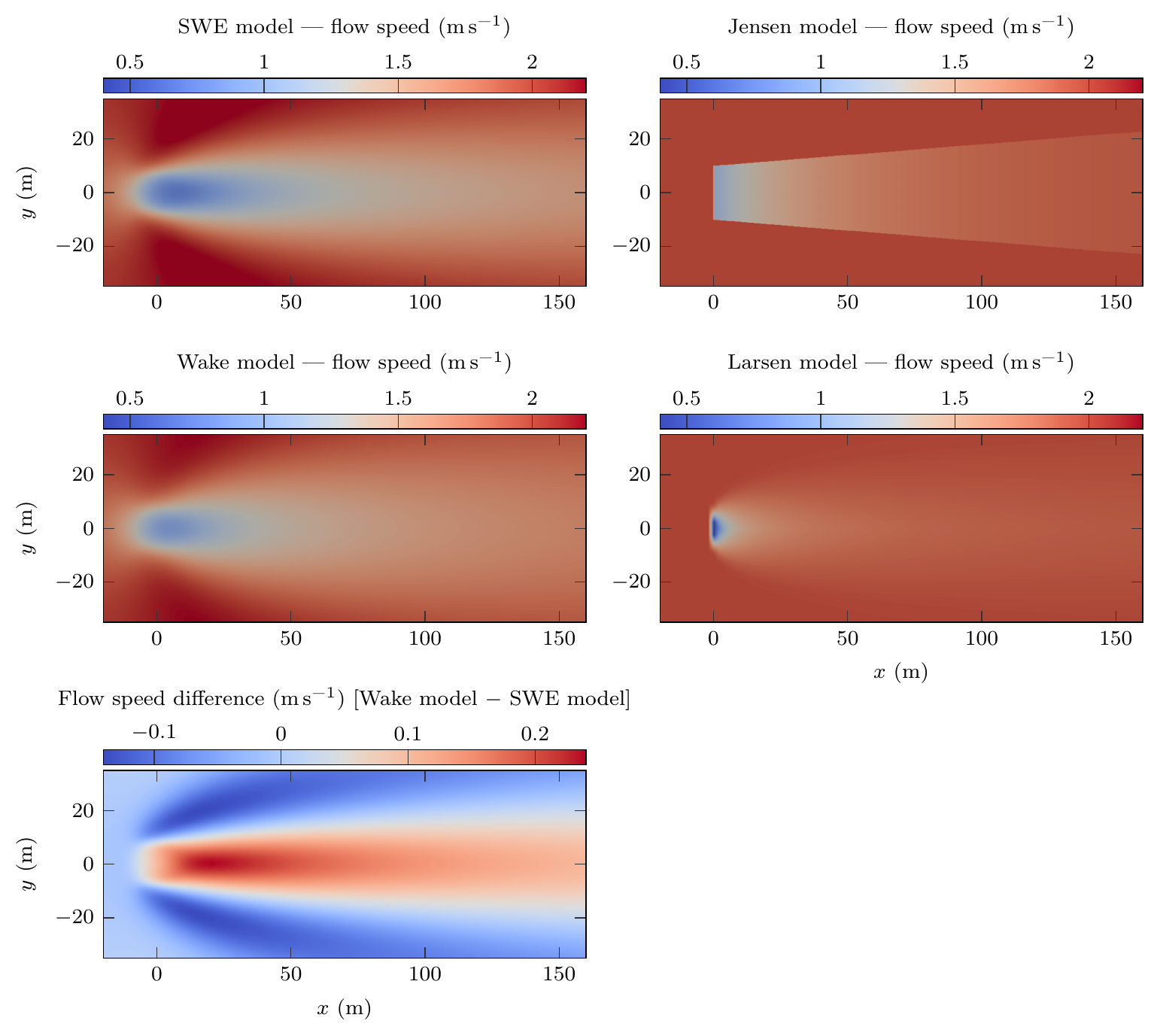}
  \caption{Wake behind a \SI{20}{\m} wide turbine placed at the origin in
    \SI{2}{\m\per\s} flow using the SWE model (top left), the wake model used
    in this work (middle left) and their difference (bottom left). Wake behind
    the same turbine is also modelled using the Jensen and Larsen models (top
    right and middle right). Each wake was modelled in a domain measuring
    \SI{1000}{\m} by \SI{400}{\m}, however only the region around the turbine
    is displayed here
  }
  \label{fig:wind-and-swe-wake}
\end{figure*}

One of the objectives of this work is to produce a simple and efficient
model to be used in global optimisation schemes, as a proxy to a
shallow water model. A natural choice for such a wake model is one
that is produced using the shallow water equations themselves. The reduction
factors required by the wake model can be generated by placing a turbine in
uniform flow, here \SI{1}{\m\per\s}, and observing the flow field around it.
The wake -- which predominantly acts in the direction parallel to the
free-stream flow with only minor components perpendicular to the free stream
flow direction -- can be used as a reduction factor ($c$ in
\cref{eq:reduced-velocity}). For simplicity, only the component parallel to
the free stream flow is extracted and taken as a reduction factor. This is
shown in \cref{fig:wind-and-swe-wake}.

\begin{figure*}[tb]
  \includegraphics[center]{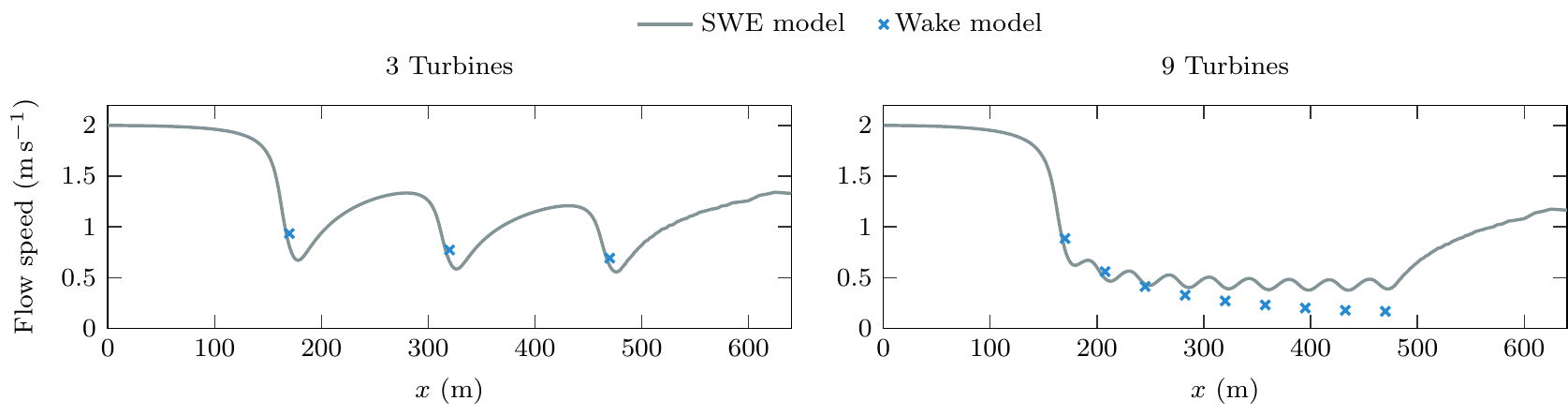}
  \caption{Combined wake for the SWE model and the
    presented wake \cref{eq:reduced-velocity} and combination
    \cref{eq:combination} models. The turbines are spaced \SI{150}{\m} and
    \SI{37.5}{\m} apart respectively for the 3 (left) and 9 (right) turbine
    cases. The inflow of \SI{2}{\m\per\s} was positioned at $x=0$ in a
    \SI{50}{\m} deep channel. Turbines were aligned parallel to the flow
    direction with the first turbine positioned at $x=170$\,\si{\m}. The error
    in combining wakes grows as more wakes are included.
  }
  \label{fig:wake-combine}
\end{figure*}

The velocity field of this approximate wake model compared to a shallow water
equation simulation with \SI{2}{\m\per\s} inflow is shown in
\cref{fig:wake-combine}. Visually, the wake appears similar. The difference
between the wakes shows that the model is largely accurate with discrepancies
lying mostly in the wake directly behind the turbine. Along the centre-line
of the wake the approximate wake model overestimates the flow whilst on either
side of the centre line it is underestimated.
Clear approximations have been made in the selection of the approximate wake
model.  However, this is deemed sufficient here as this wake model is only
used in the first, global stage of the hybrid optimisation approach, and the
optimised layout is then fed into the full shallow water model during the
second optimisation stage.  Future work could improve on this wake model by
calculating an ambient flow velocity which is representative of the turbine
site of interest before then running a shallow water equation simulation for a
single turbine in that flow velocity to generate the wake reduction factor
field or by having a `library' of models suitable for different flow speeds.

The product of the individual wake factors is taken to produce the reduction
factor
\begin{equation}
    c(x, y, \vec{m}) = \prod_{i=1}^{n} r_i(x, y, \vec{m}),
\label{eq:combination}
\end{equation}
where $c$ is the combined factor of $n$ individual reduction factors, $r_i$.
Whilst a number of combination techniques exist in the literature -- such as
the sum of squares of velocity deficits employed in the \citet{jensen1983}
model -- they often do not account for negative velocity deficits, that is
where the flow is accelerated.  The combination method \cref{eq:combination}
was chosen for its simplicity and its ability to deal with such cases. Future
work could consider an alternative wake combination methods which may be used
in cases where the flow velocity increases.

The validity of the combination model was tested by comparing the flow speed
behind turbines using the SWE and wake models. A flow of \SI{2}{\m\per\s} was
prescribed with turbines aligned in the direction of flow. For a small number
of turbines, the wake combination yields flow velocities adequately close to
that produced by the SWE model (\cref{fig:wake-combine} depicts this for 3
and 9 turbines). As the number of turbines increases, the difference between
the SWE and wake model increases. Specifically, velocities are underestimated
in the wake model. During optimisation, should a turbine be in the wake of
many others it is likely to be moved out of the wake and thus the errors
introduced in this way should be somewhat self-limiting.

\makeatletter{}\section{Optimisation schemes}\label{sec:optimisation}
\makeatletter{}Two approaches to global optimisation will be considered in this work: genetic
algorithms and a gradient-based `basin-hopping' scheme.

\makeatletter{}\subsection{Genetic algorithms}\label{subsec:genetic}

\begin{figure}[tb]
  \includegraphics[center]{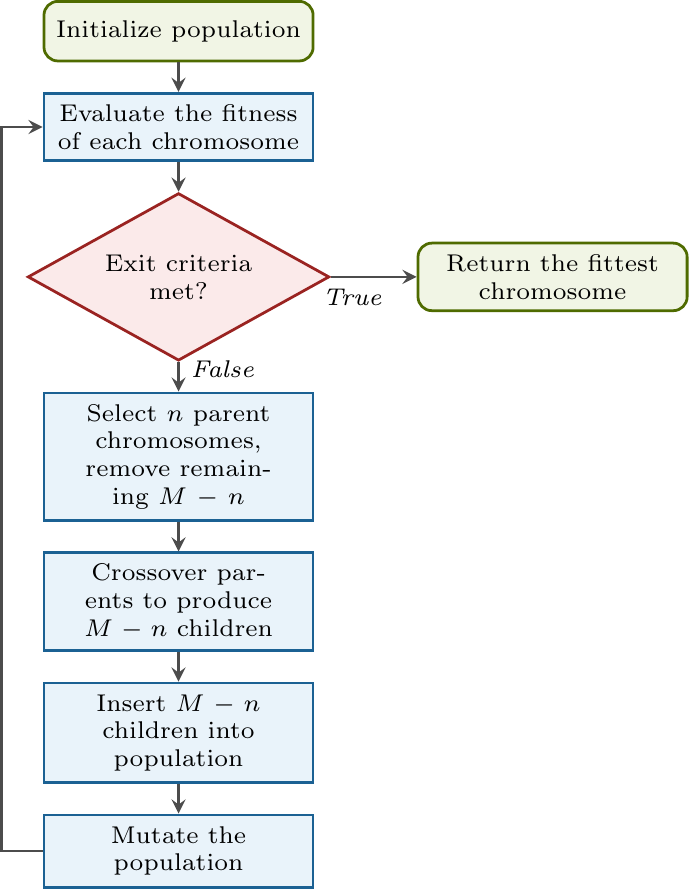}
  \caption{Typical structure of a genetic algorithm. A population of $M$
    (here, \num{100}) chromosomes is initialized, each chromosome representing
    the location of $N$ turbines. The fitness of each chromosome is evaluated
    (here the power extracted from the turbine array) so that exit criteria
    may be checked.  If the exit criteria are not met then the selection
    operator picks $n$ chromosomes (here, \num{70} based on a survival rate of
    \num{0.7}) from the population and removes the remaining $M-n$ chromosomes.
    The crossover operator takes the $n$ chosen parents and produces another
    $M-n$ children which are inserted into the population which is then
    mutated to introduce new information to the system. The cycle is repeated
    until the exit criteria are met.
  }
  \label{fig:genetic-algo}
\end{figure}

Genetic algorithms may be described as `stochastic search approaches based on
randomized operators' \citep{chelouah2000} which mimic the process of natural
selection. A genetic algorithm operates on a population containing a number of
chromosomes, each representing a vector of bits which make up the set of
variables used to evaluate the functional of interest in the optimisation
problem being posed. Here the chromosome is made up of the positions of the
turbines as per \cref{eq:m}. The fitness of a chromosome is defined to be
the evaluation of the functional of interest, that is the power production of
the farm. When applied to optimisation problems the algorithm iteratively
applies a number of genetic operators on the population until an exit criteria
is met. The typical structure of a genetic algorithm is displayed in
\cref{fig:genetic-algo}.

Three main genetic operations are performed during each iteration of the
algorithm: selection, crossover and mutation. The role of the selection
operator is to pick a number of the fittest chromosomes to propagate through
to the next generation (or iteration), the remaining chromosomes are removed
from the population. The crossover operator takes the selected chromosomes and
produces new `children' chromosomes to replace those not chosen by the
selection operator. The mutation operator acts on the new population in order
to introduce new information to the system. The new population is evaluated
and tested against exit criteria. The cycle is repeated until exit criteria
are met. Our genetic algorithm is seeded with a number of sensible
automatically generated layouts spanning different areas of the domain as well
as randomly generated layouts. The automatically generated layouts include
turbines aligned along each edge of the site, along the diagonals, and spread
across the site in arrays of varying size and dimension.  Further information
on genetic algorithms and the operators implemented and used in this work is
discussed in general in \citet{haupt2004}. Specifics regarding uniform
crossover are discussed in more detail in \citet{syswerda1989} and
\citet{chawdhry1998} whilst fitness-proportionate mutation is discussed in
\citet{pham1997}.

\makeatletter{}\subsection{Gradient-based optimisation and `basin-hopping'}\label{subsec:gradient-based}

\begin{figure}[tb]
  \includegraphics[center]{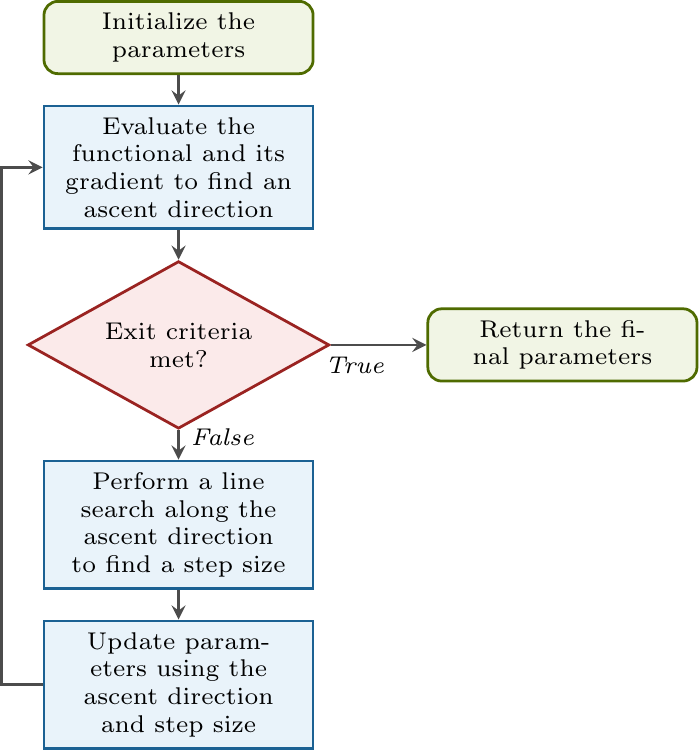}
  \caption{Typical gradient-based optimisation algorithm. The functional of
    interest to be maximised and its gradient are evaluated for input
    parameters $\vec{m}$ to find an ascent direction. The exit criteria are
    then evaluated. If they are not met then a line search is performed along
    the ascent direction to find a step size. The parameters are then updated
    accordingly using the ascent direction and step size. This cycle is
    repeated until the exit criteria are met.}
  \label{fig:gradient-algo}
\end{figure}

Genetic algorithms are global in nature and require many iterations to find an
optimised state.  This is because they depend on randomised operators to
search the full solution space. Gradient-based optimisation algorithms are
typically local in nature (see \cref{fig:schematic}) and use the gradient of
the functional of interest with respect to the parameters to guide the
algorithm's walk through the search space to an optimal solution
(\cref{fig:gradient-algo}). This can greatly decrease the number of iterations
required for convergence. Consequently this opens up the option for allowing
the use of a more costly method to evaluate the functional of interest, i.e.\
a better approximation of the real world.

\begin{figure}[tb]
  \includegraphics[center]{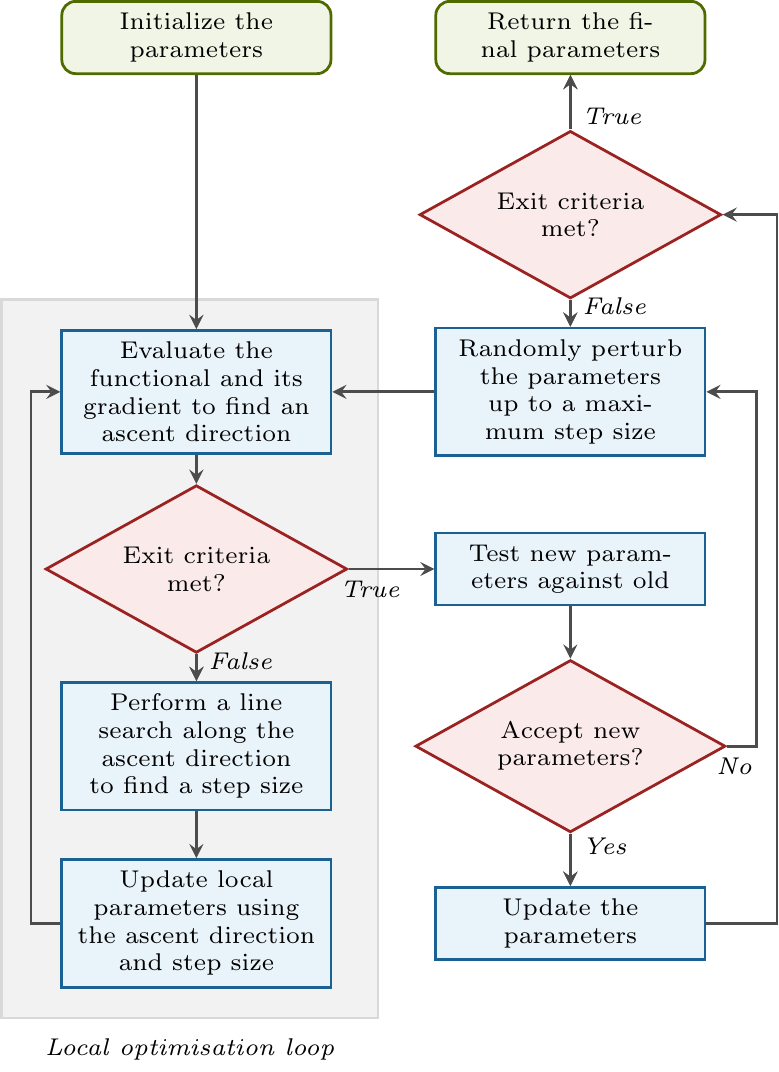}
  \caption{Schematic of the basin-hopping algorithm. After the parameters are
  initialized they are passed to a local optimization loop
  (\cref{fig:gradient-algo}). Following local optimization the new parameters
  are tested against the new parameters (using the Metropolis criterion
  \citep{li1987} of standard Monte Carlo algorithms). If the new parameters are
  accepted they are tested against the exit criteria. If the exit criteria are
  not met or the parameters were not accepted when tested they are randomly
  perturbed up to a given step size and passed back to the local optimization
  algorithm. This process is repeated until the exit criteria are met.
}
  \label{fig:basinhopping-algo}
\end{figure}

In an attempt to alleviate the limitations of a local optimisation method,
\citet{wales1997} present a stochastic `basin-hopping' algorithm which attempts
to find the global maximum of a functional by expanding on a typical
gradient-based optimisation algorithm by iteratively performing the following
steps (\cref{fig:basinhopping-algo}):
\begin{enumerate}
  \item randomly perturb the parameters,
  \item perform a local optimisation of the functional (\cref{fig:gradient-algo}),
  \item accept or reject the new parameters based on the functional value.
\end{enumerate}

The algorithm employs a variable step size (for the random perturbation) which
adjusts according to the rate at which new parameters are accepted or rejected
-- if the rate at which parameters are accepted is too high then it is likely
that the algorithm is stuck in a local maximum, and the step size is increased
in order to hop out of the maximum.  Due to the more extensive search of the
solution space, the basin-hopping algorithm often takes an order of magnitude
more iterations than local gradient-based algorithms, and this has to be taken
into account in the complexity of the underlying model used.

\makeatletter{}\subsection{Computation and validation of the functional gradient}

In this work the gradient of the total power generated by the array with
respect to the turbine coordinates is required in order to use gradient-based
optimisation. More complex functionals could of course also be used, e.g.\
those incorporating economic cost models \citep{culley2014}.  For the SWE
model the generation of the gradient information is achieved via an adjoint
calculation using \da \citep{farrell2013,funke2014}. However, this approach to
adjoint development is not applicable to the approximate wake model developed
here. Instead, automatic (or algorithmic) differentiation (AD)
\citep{griewank2008,naumann2011} is used which exploits the fact that no
matter how complex an algorithm is, it simply executes a number of simple
arithmetic operations and basic functions. AD thus repeatedly applies the
chain rule to these procedures to compute derivatives at machine precision.
This work employs the `ad' Python package \citep{pythonad} for this purpose.

Unlike a wake model represented as an explicit analytical expression (e.g.\ as
in the Jensen and Larsen models), the standard version of the AD tool cannot
be used to calculate certain components of the wake model, such as the spatial
derivatives of the wake reduction factor field or ambient flow field as they
are stored as \fenics objects which are not supported by the ad Python package.
However, thanks to the \fenics framework \citep{logg2012} upon which \otf is
built it is simple to efficiently and accurately compute these derivatives
numerically. The ad package was extended to allow access to the spatial
derivatives of the wake reduction factor calculated using \dolfin (the Python
interface to \fenics).

The implementation of the wake model gradient was verified with Taylor
remainder convergence tests, similar to \citet{funke2014}. A differentiable
functional, $J(\vec{m})$ may be expanded about a point $\vec{m}+\vec{\delta
m}$ using Taylor series:
\begin{equation}\label{eq:taylor-expansion}
  J(\vec{m}) = J(\vec{m}+\vec{\delta m}) -
  \frac{\mathrm{d}J}{\mathrm{d}\vec{m}}\delta \vec{m} +
  \mathcal{O}(\|\delta \vec{m} \|^2),
\end{equation}
thus the convergence order using first-order gradient information should
satisfy
\begin{equation}\label{eq:taylor-convergence}
  \left| J(\vec{m}+\vec{\delta m}) - J(\vec{m}) -
  \frac{\mathrm{d}J}{\mathrm{d}\vec{m}}\delta \vec{m} \right|
  = \mathcal{O}(\|\delta \vec{m} \|^2).
\end{equation}
\Cref{eq:taylor-convergence} is a strong test for checking that the
gradient computation is correctly implemented. Similarly, we expect the
convergence of the functional without gradient information to be of order
$\mathcal{O}(\|\delta \vec{m}\|^1)$.

\begin{figure}[tb]
  \includegraphics[center]{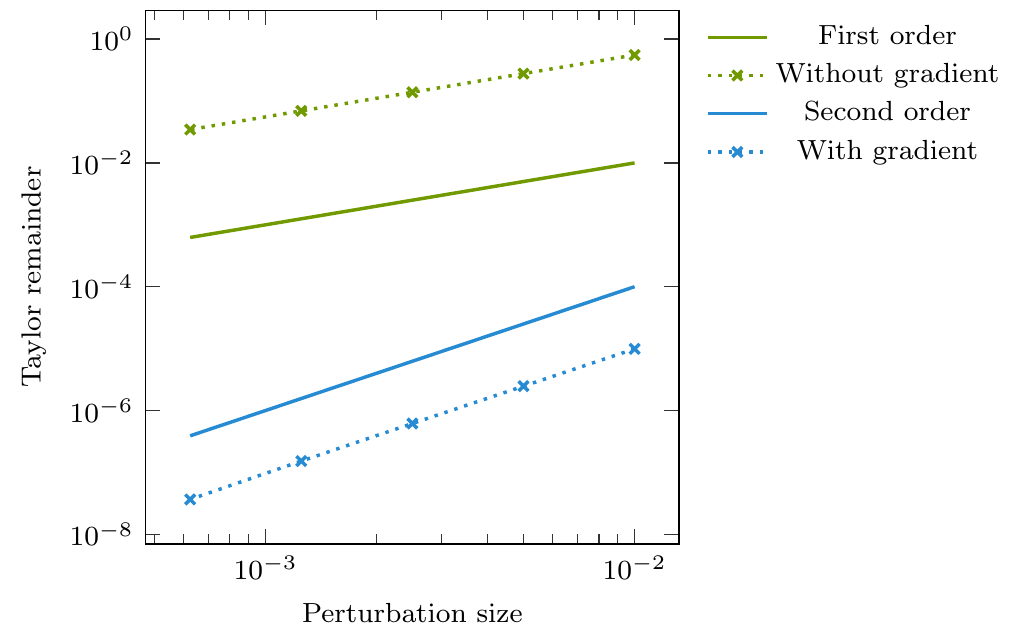}
  \caption{Expected (`First-order' and `Second-order') and achieved (`Without
    gradient' and `With gradient') orders of convergence in the Taylor
    remainder tests for the gradient of the wake model. The Taylor remainder
    is the left hand size of \cref{eq:taylor-convergence} (but only including
    the first two terms for first-order convergence) whilst the perturbation
    size is $\vec{\delta m}$ in the same equation. Turbine coordinates,
    $\vec{m}$, are perturbed in a random direction by $\vec{\delta m}$ from an
    initial set of chosen coordinates.
  }
  \label{fig:convergence}
\end{figure}

Tests were carried out for a number of flow regimes and parameter sizes
($\delta\vec{m}$) to test the validity of the gradient calculation. The
results of a test with two turbines deployed randomly in a domain measuring
\SI{640}{\m} by \SI{320}{\m} with a prescribed flow of
$(1+x/1280,1+y/640)$~\si{\m\per\s} (to ensure that the ambient flow has a
gradient) are shown in \cref{fig:convergence}. As expected, the model achieves
second-order convergence.

\makeatletter{}\section{Results and discussion}\label{sec:results}

\begin{figure*}[tb]
  \subcaptionbox{Scenario 1: a straight channel.}
      [0.30\textwidth]{\includegraphics[width=0.30\textwidth]{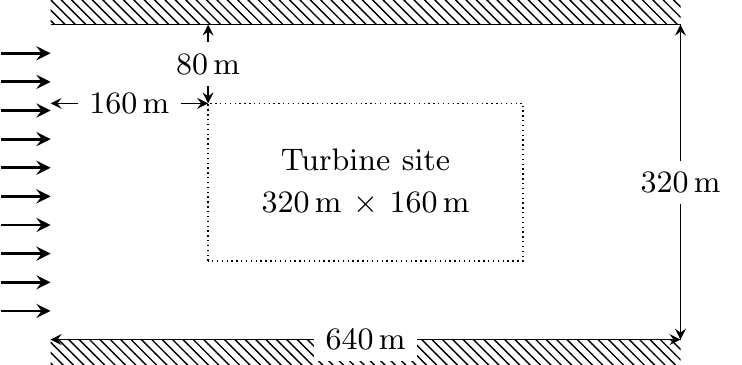}}
      \hfill
  \subcaptionbox{Scenario 3: a headland.}
      [0.30\textwidth]{\includegraphics[width=0.30\textwidth]{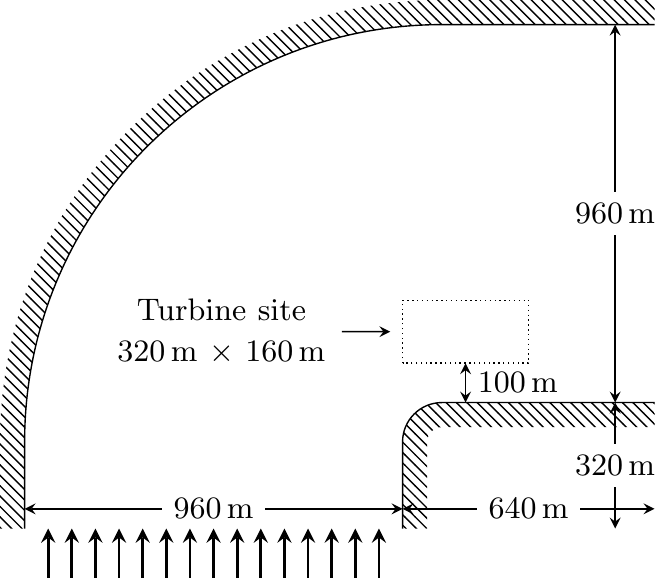}}
      \hfill
  \subcaptionbox{Scenario 4: a channel between the coast and an island.}
      [0.30\textwidth]{\includegraphics[width=0.30\textwidth]{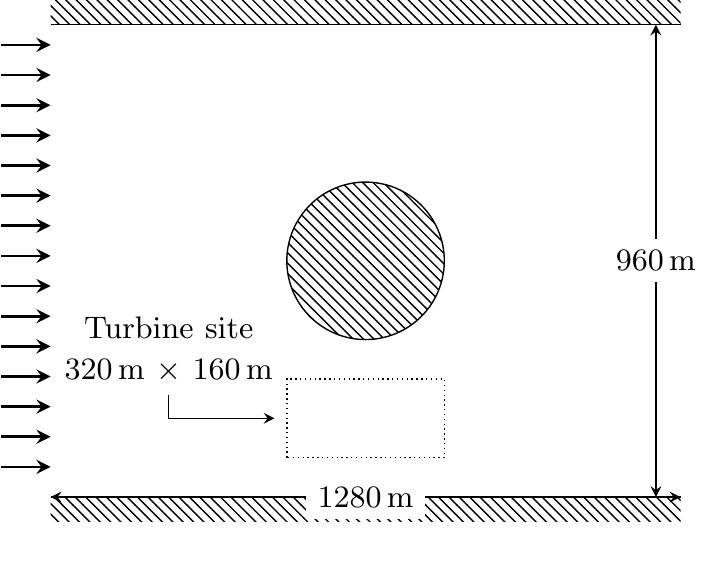}}
  \caption{Scenarios 1, 3 and 4 as per \protect\citet{funke2014} which were
    in-turn motivated by \protect\citet{draper2010}. They represent a number
    of situations which yield fast flow and thus potential tidal array sites.
    Lines of arrows indicate the inflow (a constant of \SI{2}{\m\per\s}) for
    the site whilst shaded areas indicate where there is no flow (i.e.\ land).
    A dashed box indicates this site the turbines are constrained to reside
    within.}
  \label{fig:scenarios}
\end{figure*}

\makeatletter{}Three test scenarios -- scenarios 1, 3 and 4 as per \citet{funke2014} -- were
used for testing with each turbine site measuring \SI{320}{\m} by
\SI{160}{\m}.  Schematics of these scenarios are displayed in
\cref{fig:scenarios}. In the results presented below the following shorthand
is used: `adjoint' refers to a gradient-based optimisation, using the SLSQP
algorithm \citep{kraft1988} for `local' optimisation and the basin-hopping
algorithm \citep{wales1997} for `global' optimisation. Both algorithms are
available in the \scipy Python package. Optimisation via a genetic algorithm
is denoted `genetic'. `SWE' refers to the shallow water equations model
used by \otf whilst `wake' refers to the approximate wake model developed in
this work.

Whilst local and global adjoint optimisation, and genetic optimisation may be
used with either the SWE or wake models, it is computationally prohibitive
to use global adjoint optimisation and genetic optimisation of the SWE
model as discussed in \cref{sec:introduction}. Parameters for the model and
genetic algorithm (unless otherwise stated) are displayed in
\cref{tab:parameters} and \cref{tab:genetic-parameters}.  All experiments were
run on a single core of a workstation with an Intel Xeon E3-1240 V2
(\SI{3.4}{\giga\hertz}) processor and \SI{32}{\giga B} of RAM\@.

\begin{table}[tb]
  \makeatletter{}\begin{tabular}{l l} \toprule
\textbf{Parameter}           & \textbf{Value}                      \\ \midrule
Water depth                  & $H =$ \SI{50}{\m}                   \\
Viscosity coefficient        & $\nu =$ \SI{3}{\square\m\per\s}     \\
Acceleration due to gravity  & $g =$ \SI{9.81}{\m\square\per\s}    \\
Water density                & $\rho =$ \SI{1000}{\kg\per\m\cubed} \\
Turbine friction coefficient & $K = 21$                            \\
Bottom friction coefficient  & $c_b = 0.0025$                      \\
Turbine radii                & $r =$ \SI{10}{\m}                   \\ \bottomrule
\end{tabular}

  \caption{Default parameters used in the experiments unless otherwise stated.
    Both the turbine friction coefficient and the bottom friction coefficient
    are dimensionless. The value for the turbine friction coefficient $K$ and
    the viscosity $\nu$ were motivated by \citet{funke2014}.}
  \label{tab:parameters}
\end{table}

\begin{table}[tb]
  \makeatletter{}\begin{tabular}{l l} \toprule
\textbf{Option}              & \textbf{Value}        \\ \midrule
Selection operator           & best $n$ chromosomes  \\
Crossover operator           & uniform               \\
Mutation operator            & fitness proportionate \\
Population size              & 100                   \\
Survival rate                & 0.70                  \\
Maximum mutation probability & 0.07                  \\
Maximum number of iterations & \num{10000}           \\ \bottomrule
\end{tabular}

  \caption{Parameters used for the genetic algorithm unless otherwise stated.
  The selection operator takes the best $n = 70$ chromosomes as the survival
  rate of $0.7$ is applied to a population of $100$ chromosomes. These
  parameters are largely based on trial and error and upon experience.}
  \label{tab:genetic-parameters}
\end{table}

\makeatletter{}\subsection{Timings}\label{subsec:timings}

To gain insight into the utility of the wake model, it is beneficial to have
an understanding of the relative timings for evaluations of the power
extracted from the array as well as the calculation of the gradient of the
power with respect to the turbine positions.

The main factor controlling the time to evaluate the power and its gradient in
the SWE model is the resolution of the mesh used, whilst the number of
turbines and saturation of the site (which is a function of the number of
turbines, the size of the turbines and the size of the site) determine the
time taken in the wake model. It should be noted that the ambient flow field
created for the wake model is generated by the SWE model at a cost which is
equivalent to a single power evaluation of the SWE model.  Representative
timings for the scenarios introduced above as well a more complex Pentland
Firth scenario (introduced in \cref{sec:orkney}) are presented in
\cref{tab:timing}.

Approximate timings for adjoint optimisation can be obtained by taking the
product of the number of iterations and the sum of a single power evaluation
and a single gradient evaluation\footnote{it should be noted that the number
of gradient evaluations is typically slightly lower than the number of power
evaluations in adjoint optimisation as occasionally the line search algorithm
chooses an unfavourable step length and must be  recalculated.}. For a
genetic optimisation, the approximate time taken is the product of the number
of iterations, the population size and the time for a single power evaluation.

Both the power calculation and the gradient of the wake model have an
approximate $\mathcal{O}(N^2)$ complexity as each of the $N$ turbines depends
on $N-1$ other turbines.  The wake model thus becomes expensive for large
numbers of turbines, with gradient evaluations becoming infeasible.
Evaluations of the power and gradient in the SWE model utilising the adjoint
are largely independent of the number of turbines which remains a significant
advantage for this method \citep{funke2014}. The high timings for gradient
evaluations of the wake model are due to overhead from the AD tool creating
objects for each variable associated with computing the functional of
interest. This gradient computation could be made more efficient by tailoring
the AD tool to only create objects for the variables relevant to the gradient
we are interested in.

\begin{table*}[tb]
  \begin{subtable}{\textwidth}
    \centering
    \makebox[0pt]{      \makeatletter{}\begin{tabular}{c c c c c} \toprule
\textbf{Scenario} & \textbf{Vertices} & \textbf{Power evaluation (\si{\s})} & \textbf{Gradient evaluation (\si{\s})} \\ \midrule
1                 & \num{16556}       & \num{68}                            & \num{9}                                \\
3                 & \num{31393}       & \num{167}                           & \num{19}                               \\
4                 & \num{22081}       & \num{126}                           & \num{12}                               \\
Pentland Firth    & \num{50327}       & \num{328}                           & \num{32}                               \\ \bottomrule
\end{tabular}

    }
    \caption{Representative timings for the `SWE' model utilising an
      adjoint calculation in the gradient evaluation. The power evaluation
      requires a nonlinear SWE solve and thus takes longer to compute than the
      gradient as the adjoint problem is always linear. Timings are largely
      independent of the number of turbines $N$ (here, \num{8}) but depend
      instead on the resolution of the mesh, that is the number of vertices.
    }
  \end{subtable}
  \vspace{0.5cm}

  \begin{subtable}{\textwidth}
    \centering
    \makebox[0pt]{      \makeatletter{}\begin{tabular}{c c c c} \toprule
\textbf{Turbine Site}                             & $N$       & \textbf{Power evaluation (\si{\s})} & \textbf{Gradient evaluation (\si{\s})} \\ \midrule
\multirow{3}{*}{\SI{320}{\m} by \SI{160}{\m}}     & \num{8}   & \num{0.03}                          & \num{0.6}                              \\
			                                            & \num{16}  & \num{0.09}                          & \num{4}                                \\
			                                            & \num{32}  & \num{0.3}                           & \num{25}                               \\ \addlinespace
\multirow{3}{*}{\SI{2000}{\m} by \SI{1000}{\m}}   & \num{64}  & \num{0.6}                           & \num{20}                               \\
			                                            & \num{128} & \num{2}                             & \num{155}                              \\
			                                            & \num{256} & \num{8}                             & \num{1400}                             \\ \bottomrule
\end{tabular}

    }
    \caption{Representative timings for the `wake' model. The \SI{2000}{\m} by
      \SI{1000}{\m} site represents the Pentland Firth scenario, whilst the
      \SI{320}{\m} by \SI{160}{\m} site is representative of scenarios~1, 3 and
      4.}
  \end{subtable}

  \caption{Representative timings for evaluations of power and the gradient of
    power with respect to turbine positions for the SWE model and the wake
    model. Timings for the `SWE' model are largely dependent on the mesh
    resolution whilst the size of the turbine site and
    how saturated it is with turbines heavily influence timings for the wake
    model. Turbine site saturation is a function of the number of turbines,
    the size of the turbines and the area of the turbine site. These timings
    were obtained using turbines measuring \SI{20}{\m} by \SI{20}{\m} for the
    \SI{320}{\m} by \SI{160}{\m} site and \SI{40}{\m} by \SI{40}{\m} for the
    \SI{2000}{\m} by \SI{1000}{\m} site.
  }
  \label{tab:timing}
\end{table*}

\makeatletter{}\subsection{Initial comparison}\label{subsec:initial-comparison}

\begin{figure*}[tb]
  \includegraphics[center]{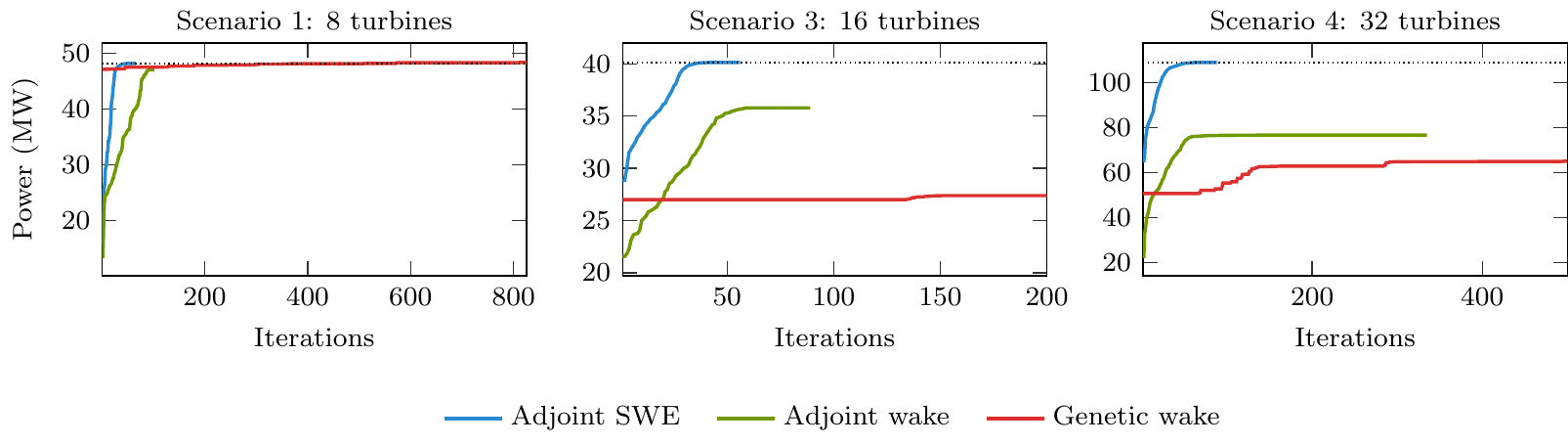}
  \caption{Power extracted from the turbine array using the SWE model
    optimised using a local gradient-based optimisation algorithm and the wake
    model using both a local-gradient based optimisation algorithm and a
    genetic algorithm.  The intermediate power was scaled using the SWE model
    in order to provide a comparable power estimate. The final power achieved
    by the SWE model is shown as a dotted line for comparison.}
  \label{fig:initial-compare}
\end{figure*}

\begin{figure*}[tb]
  \begin{subfigure}{\textwidth}
  \includegraphics[center]{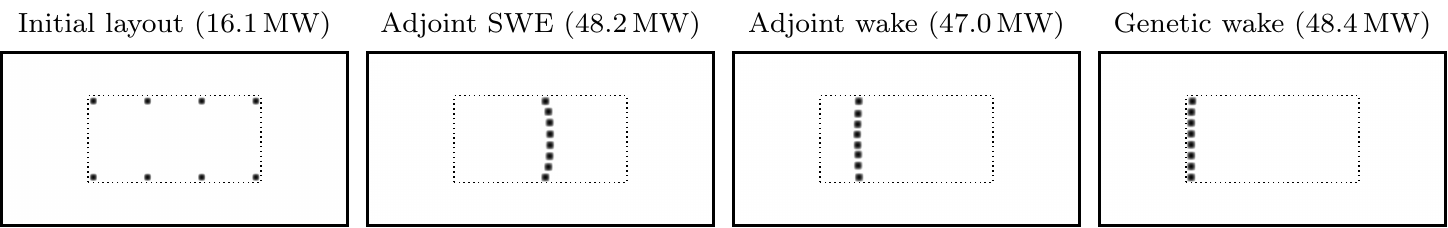}
  \caption{Optimised array layouts for 8 turbines in scenario 1. All optimised
    layouts have similar features; a line of turbines perpendicular to the
    direction of flow.
  }
  \end{subfigure}
  \vspace{1cm}

  \begin{subfigure}{\textwidth}
  \includegraphics[center]{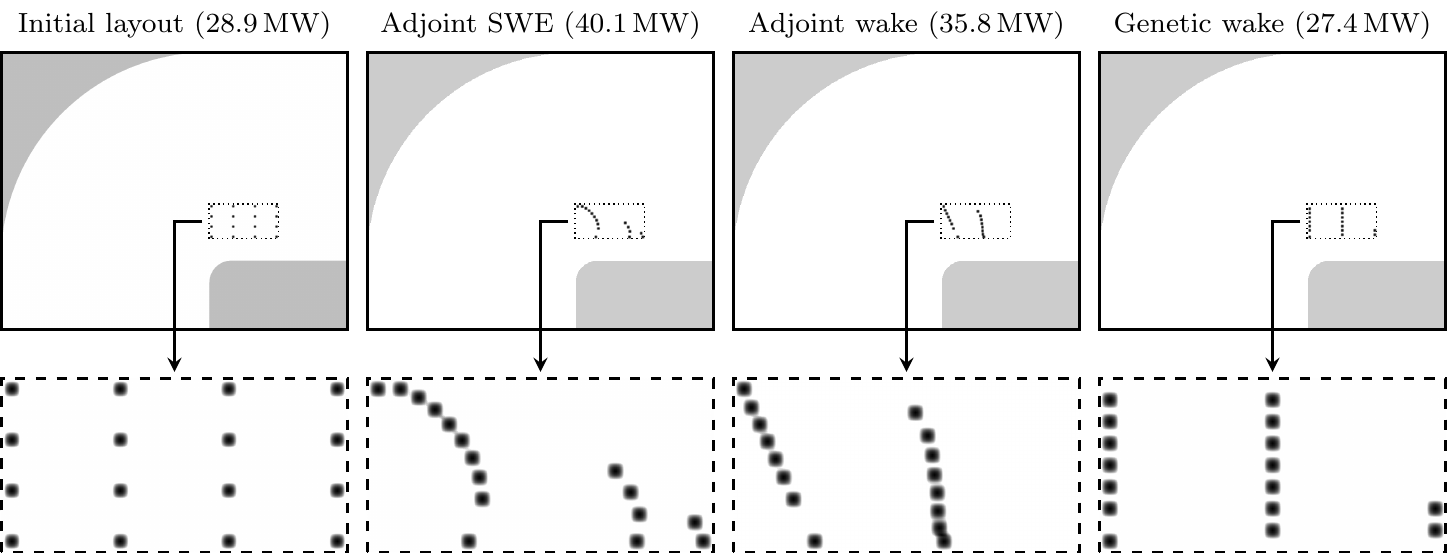}
  \caption{Optimised array layouts for 16 turbines in scenario 3. Optimised
    array layouts have similar features -- a longer line of turbines on the
    left side of the domain. These are angled approximately perpendicular to
    the direction of the flow for the adjoint SWE and adjoint wake layouts.
  }
  \end{subfigure}
  \vspace{1cm}

  \begin{subfigure}{\textwidth}
  \includegraphics[center]{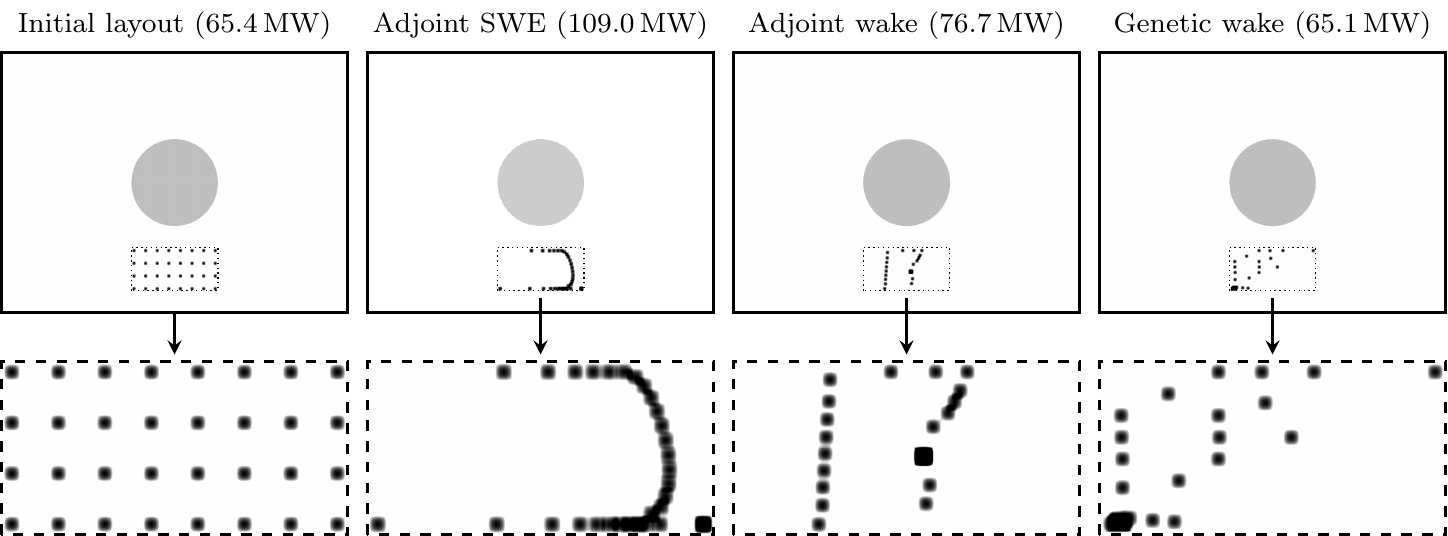}
  \caption{Optimised array layouts for 32 turbines in scenario 4. All three
    layouts contain overlapping turbines, suggesting that the site is too
    saturated with turbines. The adjoint SWE layout appears to funnel the flow
    towards the curved wall of turbines, neither of the wake model layouts
    display this property suggesting that the wake model is too simple for
    this number of turbines.
  }
  \label{fig:initial-layouts-c}
  \end{subfigure}

  \caption{Example initial and optimised layouts for three experiments from the
    initial comparison of the wake model compared to the SWE model. The
    SWE model was optimised via a local gradient-based optimisation algorithm
    (adjoint SWE). The wake model was optimised using a local gradient-based
    optimisation algorithm (adjoint wake) and a genetic algorithm (genetic
    wake).}
  \label{fig:initial-layouts}
\end{figure*}

The wake and SWE model were compared in the scenarios presented above for 8,
16 and 32 turbines (a realistic upper limit on the number of turbines that
could be deployed in a domain of this size). These initial layouts feature
turbines in regular rectangular arrrays ($4\times2$, $4\times4$, and
$8\times4$) such that turbines were spaced as far apart as the site allowed.
Results are shown in \cref{tab:initial} with layouts and optimisation
convergence plots for three examples shown in \cref{fig:initial-compare}
and~\cref{fig:initial-layouts}.

\Cref{tab:initial} shows that the SWE model consistently produces good
results: in all cases the optimised layouts produce significantly more power
than the initial layouts at iteration numbers below approximately \num{100}.
The wake model often produces good results using both the adjoint and genetic
optimisation approaches: most optimised layouts produce significantly more
power than the initial layouts. As expected, the iteration numbers are low
with the adjoint based optimisation and higher with the genetic approach.

However, the relative performance of the wake model decreases when the domain
becomes saturated with turbines. This is likely due to the simplicity of how
the wakes from multiple turbines are combined and highlights a weakness of the
wake model to approximate the SWE model when turbines are very close together.
This issue results in layouts where a number of turbines overlap. Whilst this
may yield a greater power when evaluated by the wake model the same is not
necessarily true when evaluated using the shallow water equations. This
results in certain optimised layouts being presented as extracting less power
than the initial layouts (\cref{tab:initial}, scenario~3 for \numlist{16;32}
turbines and scenario~4 for \num{32} turbines). This is evident in
\cref{fig:initial-layouts-c} where both adjoint and genetic optimisation of
the wake model yield layouts with many overlapping turbines. The extracted
power in these instances is hypothesised to be locally maximum. In some cases
turbines are also observed to overlap in the adjoint optimisation of the SWE
model.  Providing the optimisation algorithm with a minimum distance
constraint for the turbines would provide a solution to this issue but was not
done in this work.

Whilst the final power produced using the adjoint wake model is typically
lower than that obtained with the SWE model, the time taken is approximately
an order of magnitude lower. For example, in scenario 3, optimising 8 turbines
with adjoint SWE takes approximately \SI{11500}{\s} (over 3~hours) whilst the
adjoint optimisation of the wake model takes \SI{520}{\s} (under 9~minutes).

It is noted that much of the optimisation via the genetic algorithm is done in
the first few iterations and thus the number of iterations before the maximum
value is found is often relatively low (\cref{fig:initial-compare} and shown
in brackets in \cref{tab:initial}). The algorithm does not exit at this point
because the rest of the population has not converged on the same solution.
Experience indicates that this is due to being seeded with a number of
sensible layouts (similar layouts may be achieved without this seeding but in
many more iterations). Thus reducing the maximum number of iterations may
yield a very similar layout in a fraction of the time.

The results also suggest that in these flat-bottomed domains with relatively
simple ambient flow fields the solution space is simple enough to use a local
optimisation algorithm, i.e.\ basin-hopping will offer no advantage.

\begin{table*}[p]
  \makeatletter{}\begin{tabular}{c c c c c c} \toprule
\textbf{Scenario}   & $N$                 & \textbf{Optimisation} & \textbf{Iterations}     & \textbf{Power (\si{\mega\W})} \\ \midrule
\multirow{4}{*}{1} 	& \multirow{4}{*}{8} 	& Initial     	        & --                      & \num{16.08}                   \\
                   	&                    	& Adjoint SWE 	        & \num{66} 	              & \num{48.15}                   \\
                    & 			              & Adjoint wake 	        & \num{102} 	            & \num{47.04} (\num{10.20})     \\
                    &                			& Genetic wake 	        & \num{2315} (\num{45}) 	& \num{48.35} (\num{10.20})     \\ \midrule[0.1pt]
\multirow{4}{*}{1} 	& \multirow{4}{*}{16} & Initial     	        & --                      & \num{45.81}                   \\
                   	&                     & Adjoint SWE 	        & \num{88} 	              & \num{69.91}                   \\
                    & 			              & Adjoint wake 	        & \num{91} 	              & \num{60.92} (\num{12.87})     \\
                    &                			& Genetic wake 	        & \num{10000} (\num{85})  & \num{65.14} (\num{14.70})     \\ \midrule[0.1pt]
\multirow{4}{*}{1} 	& \multirow{4}{*}{32} & Initial     	        & --                      & \num{54.49}                   \\
                   	&                     & Adjoint SWE 	        & \num{102} 	            & \num{95.12}                   \\
                    &                			& Adjoint wake 	        & \num{141} 	            & \num{63.26} (\num{20.05})     \\
                    & 			              & Genetic wake 	        & \num{3009} (\num{1045}) & \num{60.26} (\num{19.70})     \\ \midrule[0.1pt]
\multirow{4}{*}{3} 	& \multirow{4}{*}{8} 	& Initial     	        & --                      & \num{20.72}                   \\
                   	&                    	& Adjoint SWE 	        & \num{63} 	              & \num{31.33}                   \\
                    & 		              	& Adjoint wake 	        & \num{85} 	              & \num{29.16} (\num{11.93})     \\
                    & 		              	& Genetic wake 	        & \num{5203} (\num{1069}) & \num{25.13} (\num{10.55})     \\ \midrule[0.1pt]
\multirow{4}{*}{3} 	& \multirow{4}{*}{16} & Initial     	        & --                      & \num{28.88}                   \\
                   	&                     & Adjoint SWE 	        & \num{56} 	              & \num{40.13}                   \\
                    & 		              	& Adjoint wake 	        & \num{89} 	              & \num{35.79} (\num{19.07})     \\
                    & 		              	& Genetic wake 	        & \num{10000}             & \num{27.39} (\num{15.08})     \\ \midrule[0.1pt]
\multirow{4}{*}{3} 	& \multirow{4}{*}{32} & Initial     	        & --                      & \num{30.95}                   \\
                   	&                     & Adjoint SWE 	        & \num{99} 	              & \num{47.37}                   \\
                    & 		              	& Adjoint wake 	        & \num{120} 	            & \num{34.80} (\num{30.33})     \\
                    & 		              	& Genetic wake 	        & \num{2738}  (\num{286}) & \num{26.98} (\num{26.69})     \\ \midrule[0.1pt]
\multirow{4}{*}{4} 	& \multirow{4}{*}{8} 	& Initial     	        & --                      & \num{30.45}                   \\
                   	&                    	& Adjoint SWE 	        & \num{52} 	              & \num{81.97}                   \\
                    & 		               	& Adjoint wake 	        & \num{82} 	              & \num{70.50} (\num{133.37})    \\
                    & 		               	& Genetic wake 	        & \num{2158} (\num{754}) 	& \num{53.10} (\num{159.84})    \\ \midrule[0.1pt]
\multirow{4}{*}{4} 	& \multirow{4}{*}{16} & Initial     	        & --                      & \num{65.00}                   \\
                   	&                     & Adjoint SWE 	        & \num{78} 	              & \num{102.88}                  \\
                    & 		              	& Adjoint wake 	        & \num{221} 	            & \num{82.94} (\num{215.39})    \\
                    & 		              	& Genetic wake 	        & \num{2459} (\num{120})  & \num{73.66} (\num{190.78})    \\ \midrule[0.1pt]
\multirow{4}{*}{4} 	& \multirow{4}{*}{32} & Initial     	        & --                      & \num{65.39}                   \\
                   	&                     & Adjoint SWE 	        & \num{88} 	              & \num{109.04}                  \\
                    & 		              	& Adjoint wake 	        & \num{335} 	            & \num{76.74} (\num{331.18})    \\
                    & 		              	& Genetic wake 	        & \num{2913} (\num{286}) 	& \num{65.10} (\num{251.41})    \\
\bottomrule \end{tabular}

  \caption{Results comparing the number of iterations and power obtained with
    the SWE and wake model for the three different scenarios considered,
    each with three numbers of turbines $N$ to be optimised. The initial
    layouts feature turbines in rectangular arrays with regular spacing making
    full use of the extent of the turbine site. Some entries for `genetic
    wake' show two iteration counts; the first represents the total number of
    iterations while the second bracketed number shows the number of
    iterations before a power within \SI{2}{\percent} of the final power was
    found. The final turbine positions obtained using the wake model were
    evaluated using the SWE model in order to provide a comparable power
    estimate in this table. The power as evaluated by the wake model is also
    shown in brackets. The particularly high power extraction in scenario~4 is
    due to the cubic dependence of power on flow speed which is close to
    \SI{4}{\m\per\s} within the turbine site due to the presence of the
    `island' in the middle of the domain.
  }
  \label{tab:initial}
\end{table*}

\makeatletter{}\subsection{Two stage global-local optimisation}\label{subsec:wake-as-guess}

Initial testing indicates that in some situations the wake model with either
genetic or adjoint optimisation may be used to get close to the power
achieved by the SWE model in a fraction of the time. The optimised layout may
be then further optimised using the more accurate SWE model. The aim of this
being that a reduction in the number of iterations required of the SWE model
would overall result in a more efficient hybrid global-local optimisation
approach whilst addressing the problem of local maxima (\cref{fig:schematic}).
The same scenarios were tested as per \cref{subsec:initial-comparison} but the
maximum number of iterations for the genetic algorithm was limited to 100 to
ensure that the initial global optimisation is cheap (for 32 turbines this
time equates to approximately 20 iterations of the SWE model).

Results are presented in \cref{tab:model2swe} and show the final optimised
power now to be more consistent across optimisation methods. There are no
distinct trends amongst the results to suggest a preferred method of
optimisation.  However, in most cases the number of iterations of the adjoint
SWE model in the hybrid approach is lower than that required by just local
optimisation of the SWE model. For these cases the total time was reduced on
average by approximately \SI{35}{\percent}. This may lead to significant
reductions in computation time when considering more realistic scenarios. A
number of results from the genetic algorithm optimisation require explanation
in light of those obtained in \cref{subsec:initial-comparison}.  Scenario~3
with \numlist{16;32} turbines and scenario~4 with 32~turbines show higher
extracted powers from \num{100}~iterations than from many thousands of
iterations as per \cref{tab:initial}. The results from \cref{tab:initial} were
highlighted above as being due to the simplicity of the wake combination model
resulting in overlapping turbines, thus when limited to \num{100}~iterations
the algorithm cannot explore the parameter space so fully, reducing the chance
of reaching such a state resulting in a greater power extraction when
evaluated by the SWE model.

In some cases the SWE model is seeded with a layout which yields a poor power
extraction (see \cref{tab:model2swe}, scenario 3 with 32 turbines, for
example) yet is able to reach a final optimised state which yields a power
similar to that achieved by the SWE model alone. This is again indicative that
these flat-bottomed domains have a solution space simple enough for a local
optimisation algorithm to be used.

\begin{table*}[tb]
  \centering
  \makebox[0pt]{    \makeatletter{}\begin{tabular}{c c c c c} \toprule
\textbf{Scenario}   & $N$                 & \textbf{Optimisation}                & \textbf{Iterations} & \textbf{Power (\si{\mega\W})}        \\ \midrule
\multirow{4}{*}{1}  & \multirow{4}{*}{8}  & Initial                              & --                  & \num{16.08}                          \\
                    &                     & Adjoint SWE                          & 66                  & \num{48.15}                          \\
                    &                     & Adjoint wake$\rightarrow$Adjoint SWE & 100, 32             & \num{47.04} (\num{10.20})$\rightarrow$\num{48.67}  \\
                    &                     & Genetic wake$\rightarrow$Adjoint SWE & 43, 26              & \num{48.99} (\num{9.96})$\rightarrow$\num{49.82}  \\ \midrule[0.1pt]
\multirow{4}{*}{1}  & \multirow{4}{*}{16} & Initial                              & --                  & \num{45.81}                          \\
                    &                     & Adjoint SWE                          & 88                  & \num{69.91}                          \\
                    &                     & Adjoint wake$\rightarrow$Adjoint SWE & 89, 55              & \num{60.92} (\num{12.87})$\rightarrow$\num{62.83}  \\
                    &                     & Genetic wake$\rightarrow$Adjoint SWE & 100, 28             & \num{67.48} (\num{14.32})$\rightarrow$\num{68.23}  \\ \midrule[0.1pt]
\multirow{4}{*}{1}  & \multirow{4}{*}{32} & Initial                              & --                  & \num{54.49}                          \\
                    &                     & Adjoint SWE                          & 102                 & \num{95.12}                          \\
                    &                     & Adjoint wake$\rightarrow$Adjoint SWE & 100, 84             & \num{64.10} (\num{19.94})$\rightarrow$\num{94.60}  \\
                    &                     & Genetic wake$\rightarrow$Adjoint SWE & 7, 35               & \num{78.49} (\num{76.47})$\rightarrow$\num{79.65}  \\ \midrule[0.1pt]
\multirow{4}{*}{3}  & \multirow{4}{*}{8}  & Initial                              & --                  & \num{20.72}                          \\
                    &                     & Adjoint SWE                          & 63                  & \num{31.33}                          \\
                    &                     & Adjoint wake$\rightarrow$Adjoint SWE & 83, 26              & \num{29.16} (\num{11.92})$\rightarrow$\num{32.34}  \\
                    &                     & Genetic wake$\rightarrow$Adjoint SWE & 81, 31              & \num{27.07} (\num{46.64})$\rightarrow$\num{29.81}  \\ \midrule[0.1pt]
\multirow{4}{*}{3}  & \multirow{4}{*}{16} & Initial                              & --                  & \num{28.88}                          \\
                    &                     & Adjoint SWE                          & 56                  & \num{40.13}                          \\
                    &                     & Adjoint wake$\rightarrow$Adjoint SWE & 87, 48              & \num{35.79} (\num{19.07})$\rightarrow$\num{44.16}  \\
                    &                     & Genetic wake$\rightarrow$Adjoint SWE & 100, 46             & \num{29.81} (\num{74.54})$\rightarrow$\num{43.77}  \\ \midrule[0.1pt]
\multirow{4}{*}{3}  & \multirow{4}{*}{32} & Initial                              & --                  & \num{30.95}                          \\
                    &                     & Adjoint SWE                          & 99                  & \num{47.37}                          \\
                    &                     & Adjoint wake$\rightarrow$Adjoint SWE & 100, 102            & \num{34.80} (\num{30.30})$\rightarrow$\num{50.47}  \\
                    &                     & Genetic wake$\rightarrow$Adjoint SWE & 100, 135            & \num{32.76} (\num{111.07})$\rightarrow$\num{50.86}  \\ \midrule[0.1pt]
\multirow{4}{*}{4}  & \multirow{4}{*}{8}  & Initial                              & --                  & \num{30.45}                          \\
                    &                     & Adjoint SWE                          & 52                  & \num{81.97}                          \\
                    &                     & Adjoint wake$\rightarrow$Adjoint SWE & 80, 80              & \num{70.50} (\num{26.54})$\rightarrow$\num{81.88}  \\
                    &                     & Genetic wake$\rightarrow$Adjoint SWE & 58, 24              & \num{80.59} (\num{145.22})$\rightarrow$\num{81.91}  \\ \midrule[0.1pt]
\multirow{4}{*}{4}  & \multirow{4}{*}{16} & Initial                              & --                  & \num{65.00}                          \\
                    &                     & Adjoint SWE                          & 78                  & \num{102.88}                         \\
                    &                     & Adjoint wake$\rightarrow$Adjoint SWE & 100, 57             & \num{82.62} (\num{43.13})$\rightarrow$\num{102.78} \\
                    &                     & Genetic wake$\rightarrow$Adjoint SWE & 100, 74             & \num{67.31} (\num{153.69})$\rightarrow$\num{104.30} \\ \midrule[0.1pt]
\multirow{4}{*}{4}  & \multirow{4}{*}{32} & Initial                              & --                  & \num{65.39}                          \\
                    &                     & Adjoint SWE                          & 88                  & \num{109.04}                         \\
                    &                     & Adjoint wake$\rightarrow$Adjoint SWE & 100, 101            & \num{76.71} (\num{64.18})$\rightarrow$\num{109.48} \\
                    &                     & Genetic wake$\rightarrow$Adjoint SWE & 100, 97             & \num{73.11} (\num{206.38})$\rightarrow$\num{107.13} \\
\bottomrule \end{tabular}

  }
  \caption{Results comparing two-stage optimisation of the wake model to a
  single local optimisation of the SWE model. Two-stage optimisation consists
  of global optimisation of the wake model (adjoint or genetic) followed by
  local optimisation of the SWE model. This it indicated using the arrow
  notation ($\rightarrow$).
}
  \label{tab:model2swe}
\end{table*}

\makeatletter{}\subsection{Inclusion of bathymetry}\label{subsec:bathymetry}

The importance of global approaches to optimisation increases with the
complexity of the solution space. In tidal turbine array optimisation one
important way this occurs is through the addition of non-flat bathymetry. The
same optimisation strategies as in \cref{subsec:wake-as-guess} were used.
However, in this section the basin-hopping algorithm (which makes use of the
local gradient-based approach) was used in the adjoint optimisation of the
wake model. The number of iterations allowed by the genetic algorithm was also
increased from \num{100} to \num{1000} to account for the more complex
solution space.

\begin{figure*}[tb]
  \includegraphics[center]{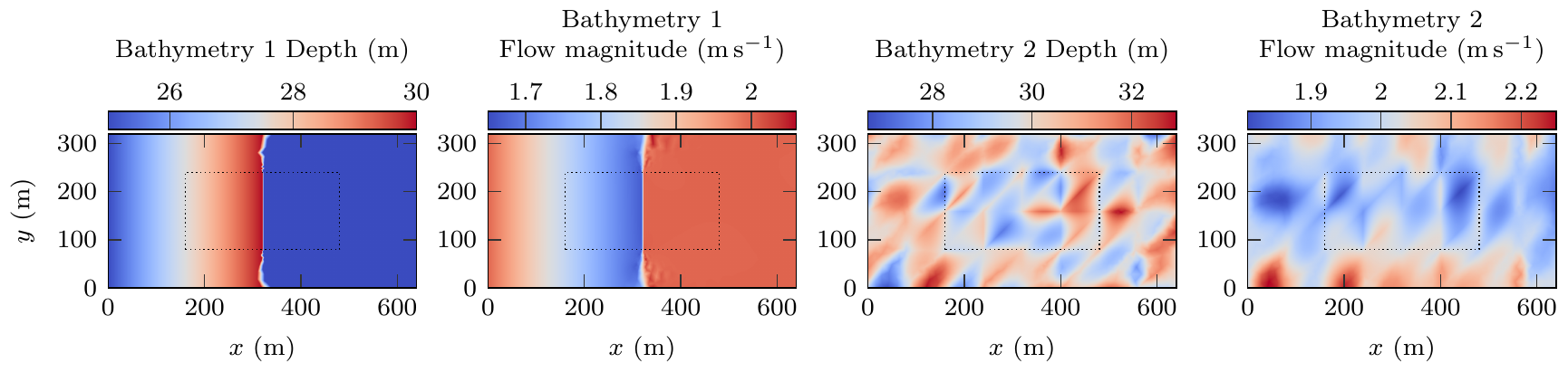}
  \caption{Bathymetry and flow magnitude due to a \SI{2}{\m\per\s} flow from the
    left of the domain for bathymetry~1 and~2. The turbine site is marked by
    the dotted rectangle.}
  \label{fig:bath12}
\end{figure*}

Three bathymetry cases were considered. The first bathymetry field is applied
to the scenario 1 channel domain from \cref{sec:results} and features depth
increasing linearly from \SIrange{25}{30}{\m} from the inflow boundary to
midway through the domain, at this point the depth decreases rapidly to
\SI{25}{\m} where it remains constant to the outflow. Thus the fastest current
-- and hence best location for the turbines -- is on the right of the domain,
where it is shallowest. However, there is a decrease in flow velocity in the
first half of the domain as one moves away from the inflow and thus turbines
in this region will not be able to reach the optimal part of the domain using
a local algorithm. The bathymetry and ambient flow are displayed in
\cref{fig:bath12}.

\begin{figure*}[tb]
  \includegraphics[center]{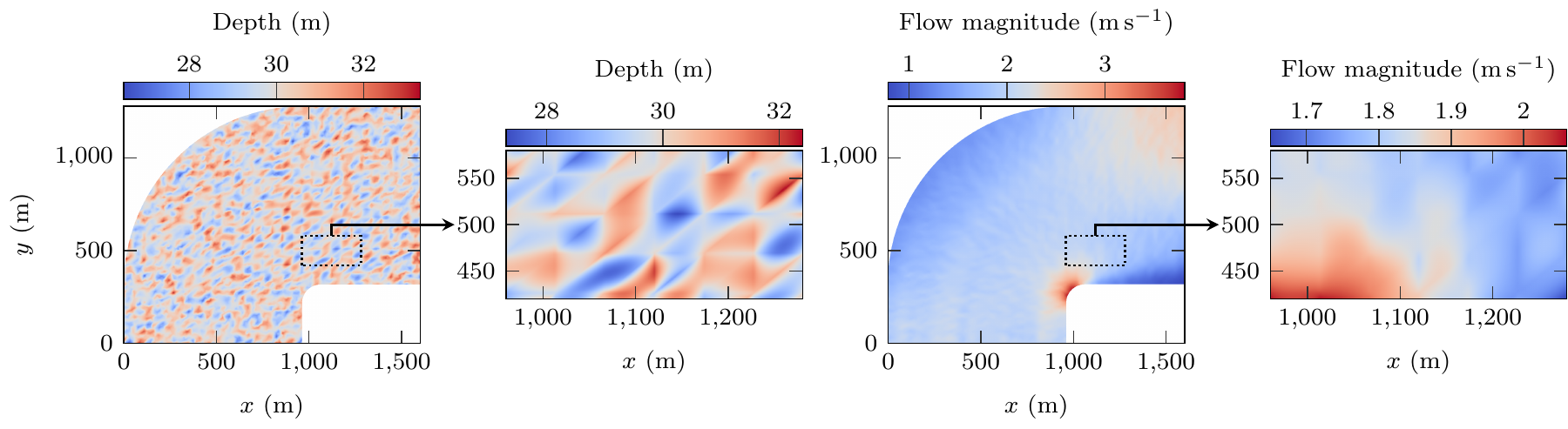}
  \caption{Bathymetry and flow magnitude due to a \SI{2}{\m\per\s} flow from the
    bottom edge of the domain for the whole domain (top) and turbine site
    (bottom) of bathymetry~3. The turbine site is marked by the dotted
    rectangle.}
  \label{fig:bath3}
\end{figure*}

The second bathymetry, again applied to the scenario 1 domain, contains
a randomly generated bathymetry. Unlike bathymetry 1, finding the optimised
layout is non-trivial and must be guided by more than intuition.  Bathymetry 3
has similar features to bathymetry 2 but is applied to a different domain
(scenario 3 as per \cref{sec:results}). The resulting flow field and
bathymetry are presented in \cref{fig:bath3}. The flow speed in the turbine
site appears more heavily influenced by the geometry of the domain than the
bathymetry, resulting in a simpler solution space than bathymetry~1 and~2.

\begin{table*}[tb]
  \makebox[0pt]{    \makeatletter{}\newcommand{\raas}{$\rightarrow$Adjoint SWE}
\newcommand{\ra}{$\rightarrow$}
\begin{tabular}{c c c c c} \toprule
\textbf{Bathymetry} & $N$                & \textbf{Optimisation} & \textbf{Iterations}  & \textbf{Power (\si{\mega\W})} \\ \midrule
\multirow{4}{*}{1}  & \multirow{4}{*}{4} & Initial               & --                   & \num{5.69}                    \\
                    &                    & Adjoint SWE           & \num{46}             & \num{7.97}                    \\
                    &                    & Adjoint wake\raas     & \num{1873}, \num{45} & \num{8.87} (\num{15.45})\ra\num{9.12}       \\
                    &                    & Genetic wake\raas     & \num{448}, \num{23}  & \num{9.09} (\num{19.16})\ra\num{9.18}       \\ \midrule[0.1pt]
\multirow{4}{*}{2}  & \multirow{4}{*}{4} & Initial               & --                   & \num{7.05}                    \\
                    &                    & Adjoint SWE           & \num{7}              & \num{8.56}                    \\
                    &                    & Adjoint wake\raas     & \num{1919}, \num{23} & \num{10.15} (\num{19.36})\ra\num{10.52}     \\
                    &                    & Genetic wake\raas     & \num{898}, \num{35}  & \num{10.13} (\num{19.94})\ra\num{10.72}     \\ \midrule[0.1pt]
\multirow{4}{*}{3}  & \multirow{4}{*}{4} & Initial               & --                   & \num{6.79}                    \\
                    &                    & Adjoint SWE           & \num{48}             & \num{9.58}                    \\
                    &                    & Adjoint wake\raas     & \num{1915}, \num{20} & \num{8.50} (\num{24.50})\ra\num{9.63}       \\
                    &                    & Genetic wake\raas     & \num{614}, \num{19}  & \num{7.80} (\num{22.57})\ra\num{9.80}       \\
\bottomrule \end{tabular}

  }
  \caption{Results comparing two-stage optimisation of the wake model to local
    optimisation of the SWE model for three cases including the effects of
    bathymetry. Complex bathymetry is likely to cause a number of local
    maxima, thus global optimisation is essential to finding the optimal
    layout.  Two-stage optimisation consists of global optimisation of the
    wake model (basin-hopping or genetic) followed by local optimisation of
    the SWE model. This is indicated using the arrow notation ($\rightarrow$).
    In bathymetry~2 the initial layout is hypothesised to be very close to a
    local maximum which is found in seven iterations by the adjoint SWE
    optimisation. For the two hybrid optimisations of the same bathymetry the
    number of adjoint SWE iterations is greater as the layout provided by the
    global optimisation of the wake model is further from a maximum.
  }
  \label{tab:bathymetry}
\end{table*}

\begin{figure*}[tb]
  \includegraphics[center]{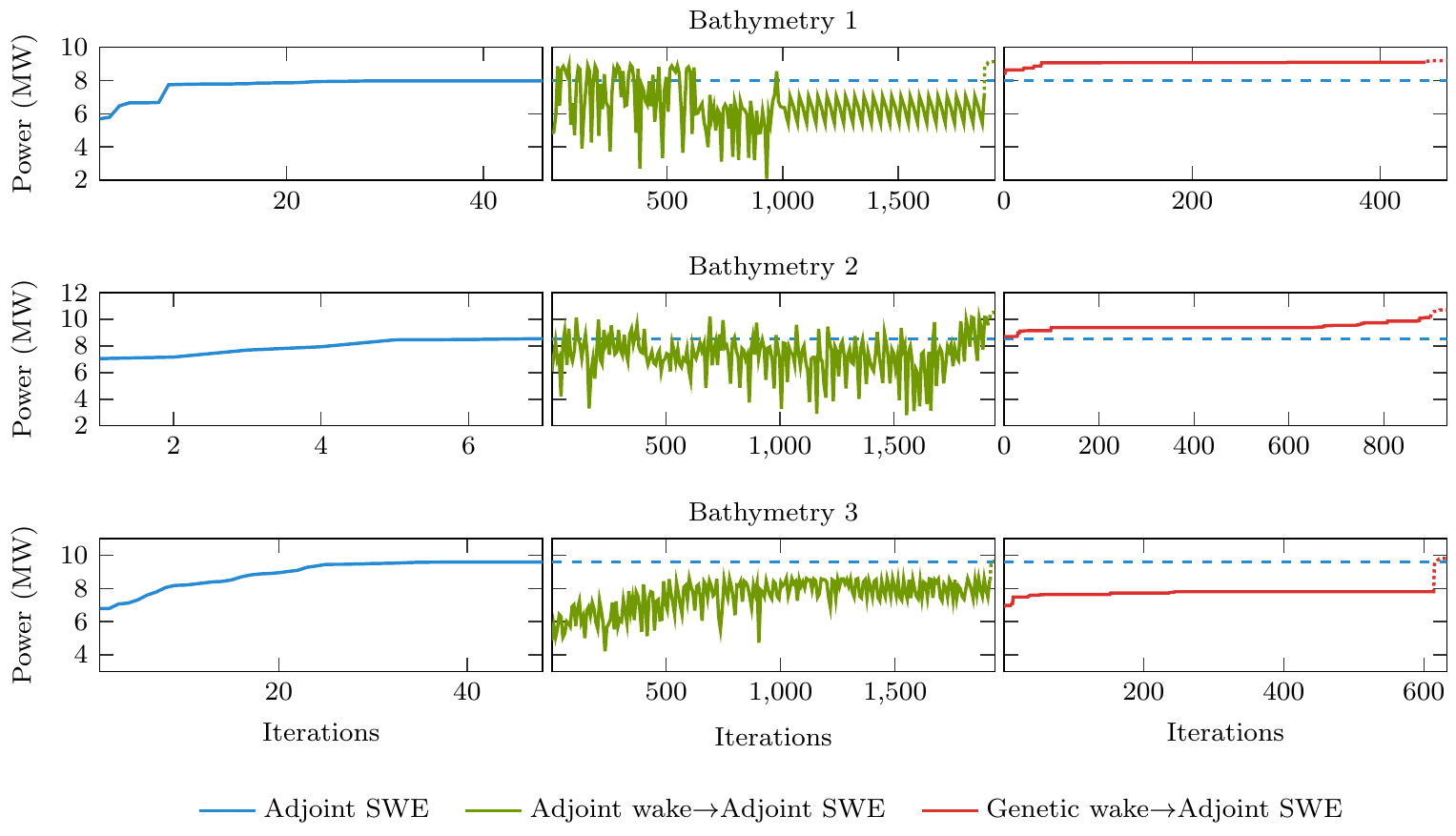}
  \caption{Optimisation convergence for bathymetry 1, 2 and~3. The second
    stage of the optimisation is dotted. A horizontal dashed line representing
    the final power of the adjoint SWE optimisation is also shown for
    comparison.  The first stages of the adjoint and
    genetic wake optimisation (shown as solid lines) are scaled such that the
    final power of the first stage is equal to the first power achieved in the
    second stage. Thus the power values displayed for the first stage merely
    show the relative improvement of the optimisation of the wake model and
    not the actual power extracted.  During the first stage of adjoint
    optimisation in bathymetry~1 after approximately \num{1000} iterations the
    basin-hopping algorithm is unable to escape a local maxima; this is due to
    the bathymetry of the domain. In this case all turbines are in the right
    part of the domain where the depth is constant and the turbines are
    perturbed such that they have little effect on each other. The optimised
    state occurs near to \num{600} iterations where the turbines are in the
    same part of the domain but are working together to increase the flow
    speed for the central turbines.
  }
  \label{fig:bath-power}
\end{figure*}

\begin{figure*}[tb]
  \begin{subfigure}{\textwidth}
  \includegraphics[center]{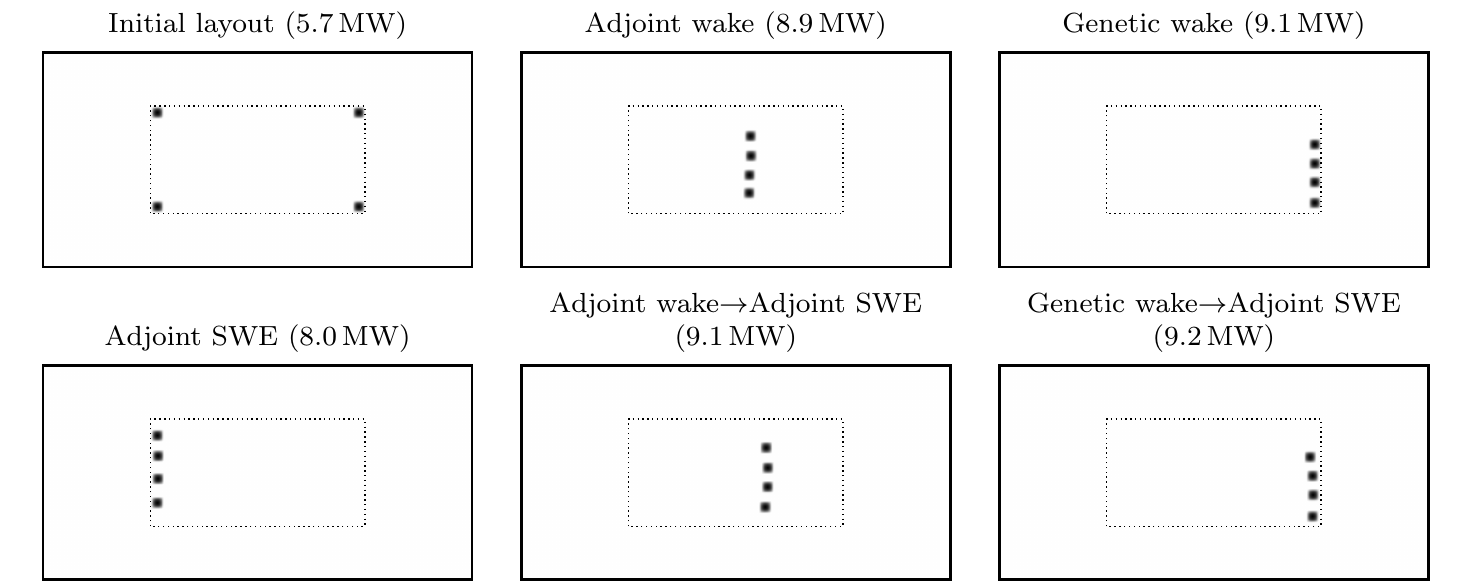}
  \caption{Optimised array layouts for 4 turbines in bathymetry~1.}
  \label{fig:bath1-layouts}
  \end{subfigure}
  \vspace{1cm}

  \begin{subfigure}{\textwidth}
  \includegraphics[center]{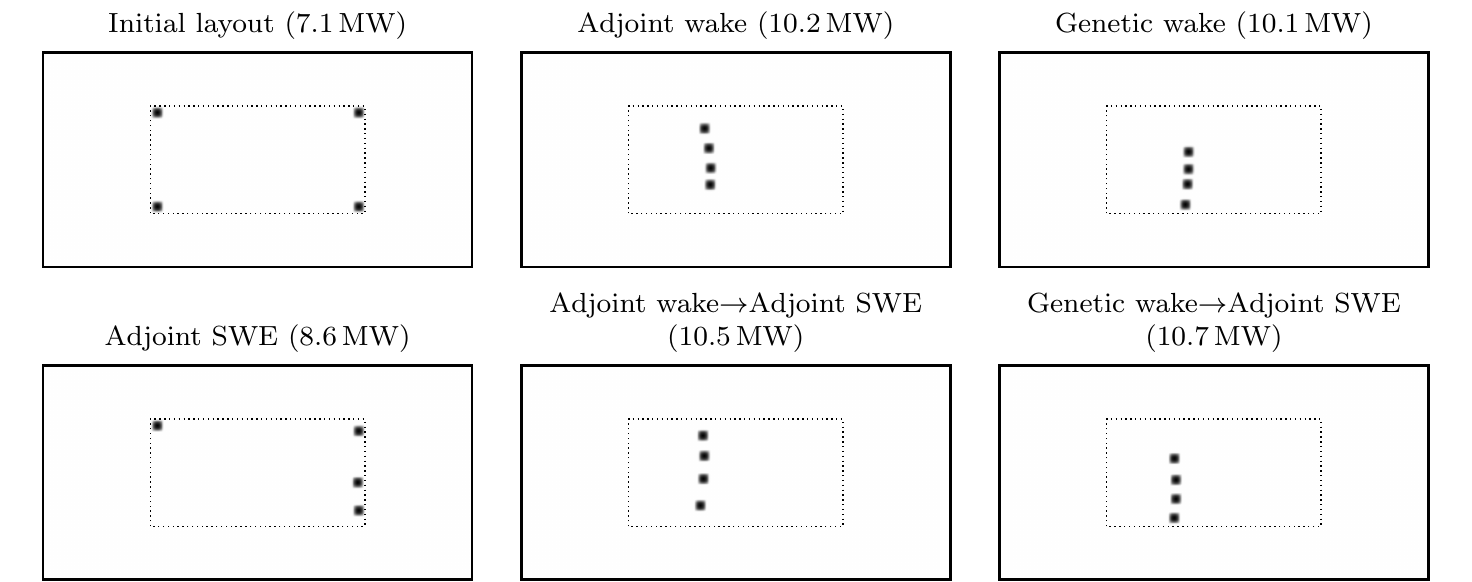}
  \caption{Optimised array layouts for 4 turbines in bathymetry~2.}
  \label{fig:bath2-layouts}
  \end{subfigure}
  \caption{Initial and optimised layouts for bathymetry scenarios~1 and~2.}
  \label{fig:bath1-2-layouts}
\end{figure*}

\begin{figure*}
  \includegraphics[center]{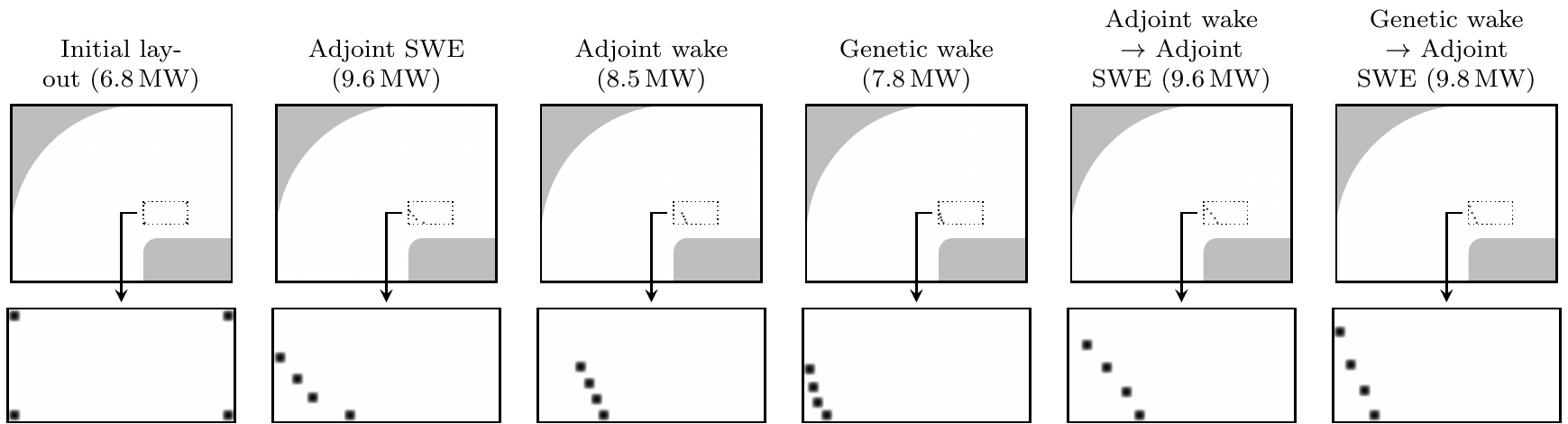}
  \caption{Optimised array layouts for 4 turbines in bathymetry~3.}
  \label{fig:bath3-layouts}
\end{figure*}

\cref{tab:bathymetry} displays the results of comparing the SWE model and the
wake model in the three different bathymetry scenarios presented above.
Convergence plots are presented in \cref{fig:bath-power} whilst initial and
optimised layouts are displayed in
\cref{fig:bath1-2-layouts,fig:bath3-layouts}.  Bathymetry~1 and 2 show that
using global optimisation schemes (basin-hopping and genetic) can yield
significant power improvements over local optimisation schemes whilst a
two-stage optimisation yields further improvements. For bathymetry 1, layout
optimisation of the SWE model (using the adjoint approach without
basin-hopping) results in all turbines on the left hand side of the domain
where the flow is slowest whilst both optimisations of the wake model (genetic
and using the adjoint approach with basin-hopping) result in turbines on the
right hand side of the domain in the faster flowing water
(\cref{fig:bath1-layouts}).  After the second stage of optimisation
\SI{14.4}{\percent} and \SI{15.2}{\percent} improvements are yielded over
local optimisation of the SWE model for the adjoint and genetic approaches
respectively.

In bathymetry~2 the local optimisation of the SWE model yields a layout
somewhat similar to the initial layout (\cref{fig:bath2-layouts}); we
hypothesise that the algorithm is stuck in a local maximum. Both global
optimisations of the wake model (genetic and basin-hopping) result in similar
layouts. After second stage optimisation the basin-hopping and genetic
approaches yield \SI{22.9}{\percent} and \SI{25.2}{\percent} improvements over
local optimisation of the SWE model.

For bathymetry~3 similar power outputs are achieved for all optimisation
approaches. The power achieved by the wake model is \SI{0.5}{\percent} and
\SI{2.3}{\percent} greater than optimisation of the SWE model for the adjoint
and genetic approaches respectively. Importantly, the layouts discovered using
the two-stage approach are achieved in approximately \SI{40}{\percent} less
time than those found by the adjoint SWE optimisation (excepting bathymetry~2
where it is hypothesised that the initial layout is very close to a local
maximum which is found in very few iterations by the adjoint SWE
optimisation).

\makeatletter{}\section{Array optimisation in the Inner Sound of Pentland Firth}\label{sec:orkney}

The presence of the Orkney Islands accelerates tidal flow through the Pentland
Firth separating the Islands from the north east coast of Scotland.  This is
therefore a site of major interest for tidal turbine array development. The
Inner Sound of the Pentland Firth lies in the channel between Stroma Island
and Caithness where the water is shallower than the surrounding areas which
further accelerates the flow and provides a potential site of appropriate depth
for turbine deployment.  It is one of the most promising locations for tidal
power in the UK and is currently under development by MeyGen Ltd.\ who have
recently been given permission to deploy a first stage of turbines expected to
produce $\sim$\SI{86}{\mega\watt}, with a long term goal of a
$\sim$\SI{400}{\mega\watt} array.

This site will be used here as a more realistic scenario to demonstrate the
work presented above.  Layout optimisation of 64 turbines (which is
approximately the number required to deliver \SI{86}{\mega\watt} of power)
will be compared for the local and hybrid global-local approaches as per
\cref{subsec:bathymetry}.

A number of parameters from \cref{tab:parameters} were adjusted to reduce
computation time for this larger, more complex problem than was previously
considered. The viscosity coefficient, $\nu$ was increased to
\SI{100}{\m\per\s\squared} and the turbine radii to \SI{20}{\m}. Depth
throughout the domain was set using a mixture of the best bathymetry data
available in different parts of the domain: from GEBCO \citep{becker2009},
Digimap \citep{digimap} and the \citet{scotlandgov}.  Accurate shoreline data
comes from the GSHHS database \citep{wessel1996}. The mesh was created using
Gmsh \citep{geuzaine2009} and consists of $\sim$\num{1.0e5} triangles varying from
\SIrange{8.16}{223.14}{\m} in size. The domain is the same as used by
\citet{funke2014} but includes bathymetry and employs a coarser mesh to reduce
computation time. Only the idealised steady flood tide of \SI{2}{\m\per\s}
along the left boundary of the domain is considered in this optimisation.

\begin{figure*}[tb]
  \includegraphics[center]{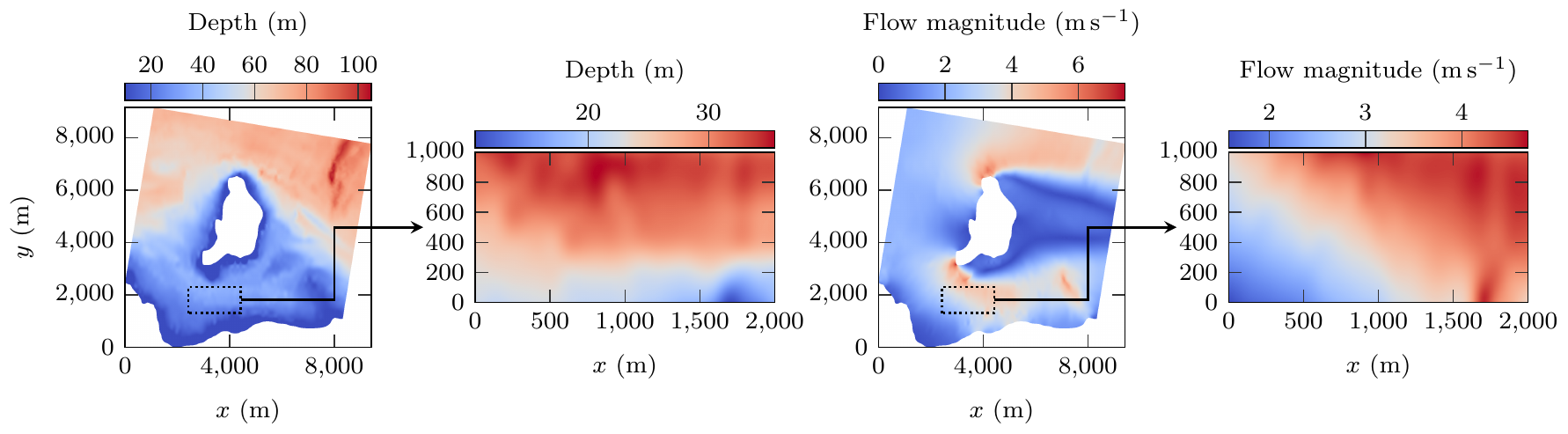}
  \caption{Bathymetry and flow magnitude for the whole domain (top) and the
    turbine site (bottom) for the Pentland Firth scenario. The turbine site is
    shown by the dotted rectangle in the top images.}
  \label{fig:orkney-bathymetry}
\end{figure*}

\begin{figure}[tb]
  \includegraphics[center]{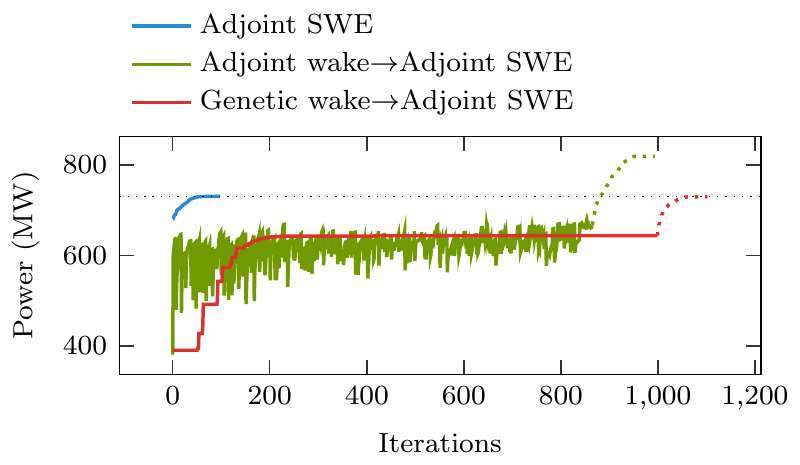}
  \caption{Optimisation convergence for 64 turbines in the Pentland Firth
  scenario. The second stages of optimisation are shown as dotted lines. The
  two stage adjoint optimisation (adjoint optimisation of the wake model
  followed by adjoint optimisation of the SWE model) yields a significant
  improvement over the adjoint SWE optimisation. The first stages of the
  adjoint and genetic wake optimisation (solid lines) are scaled such that the
  final value equals that of the same layout evaluated using the SWE
  model. Extracted powers for the first stage thus indicate the relative
  improvement during wake model optimisation and do not necessarily relate to
  the power which would be extracted from the same layout evaluated by the SWE
  model.}
\label{fig:orkney-power}
\end{figure}

The bathymetry and resulting ambient flow field for the whole domain and the
turbine site are displayed in \cref{fig:orkney-bathymetry}, along with the
allowed turbine site which is simplified compared to that being considered by
MeyGen Ltd. The power extracted from the array configurations is displayed in
\cref{tab:orkney} with optimisation convergence displayed in
\cref{fig:orkney-power}. The adjoint two-stage approach yields a 12.2\%
improvement in power over local optimisation of the SWE model alone, whilst
the two-stage genetic approach results in a power similar to local
optimisation of the SWE model. It should be emphasised that there are still
many simplifications present in the turbine and flow representations and the
figures presented need to be viewed in this light.

\begin{figure*}[tb]
  \includegraphics[center]{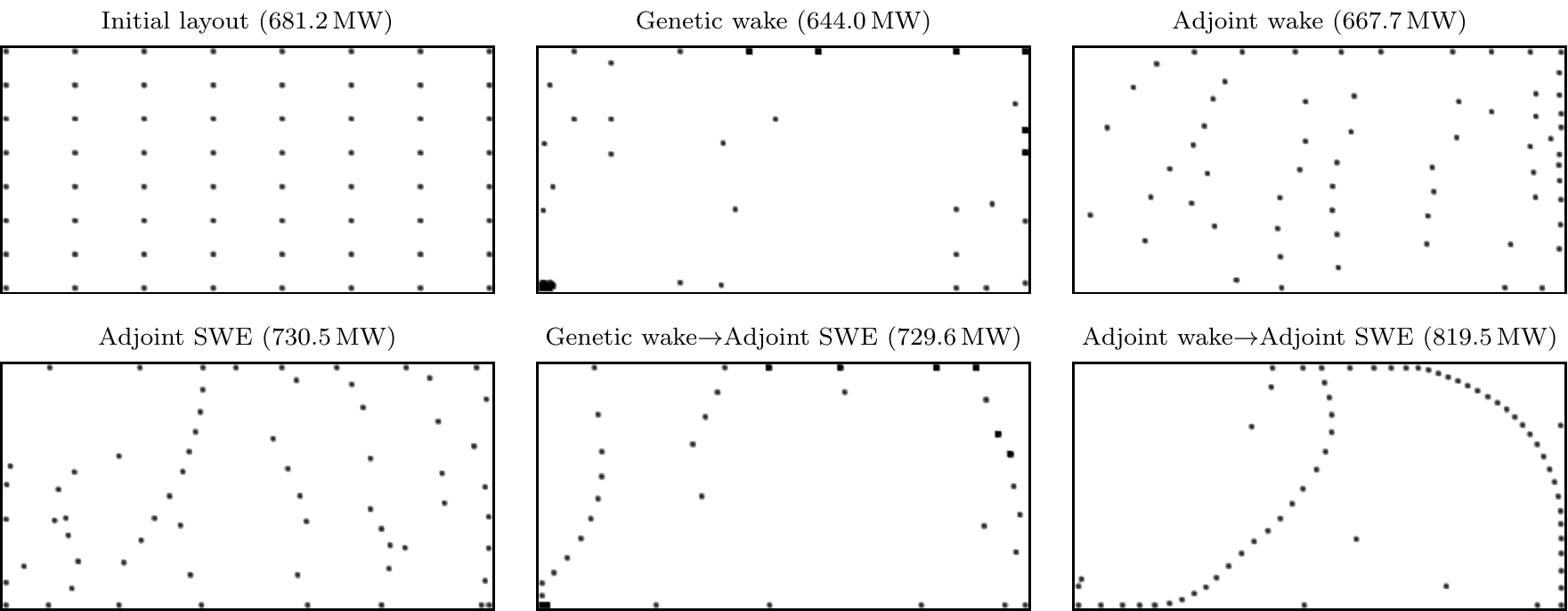}
  \caption{Initial and optimised layouts for the Pentland Firth scenario. The
    adjoint~wake$\rightarrow$adjoint~SWE optimisation represents the optimal
    array design which contains a number of barrage-like structures to channel
    the flow through the turbine array. It is noted that the genetic wake
    layout results in a number of turbines on top of each other. There are a
    number of reasons which may cause this: inadequacy of the wake model to
    approximate the SWE model when turbines are in close proximity and the
    poor suitability of genetic algorithms when optimising for many turbines
    as the parameter space grows greatly.
}
  \label{fig:orkney-layouts}
\end{figure*}

\begin{figure*}[tb]
  \includegraphics[width=\textwidth]{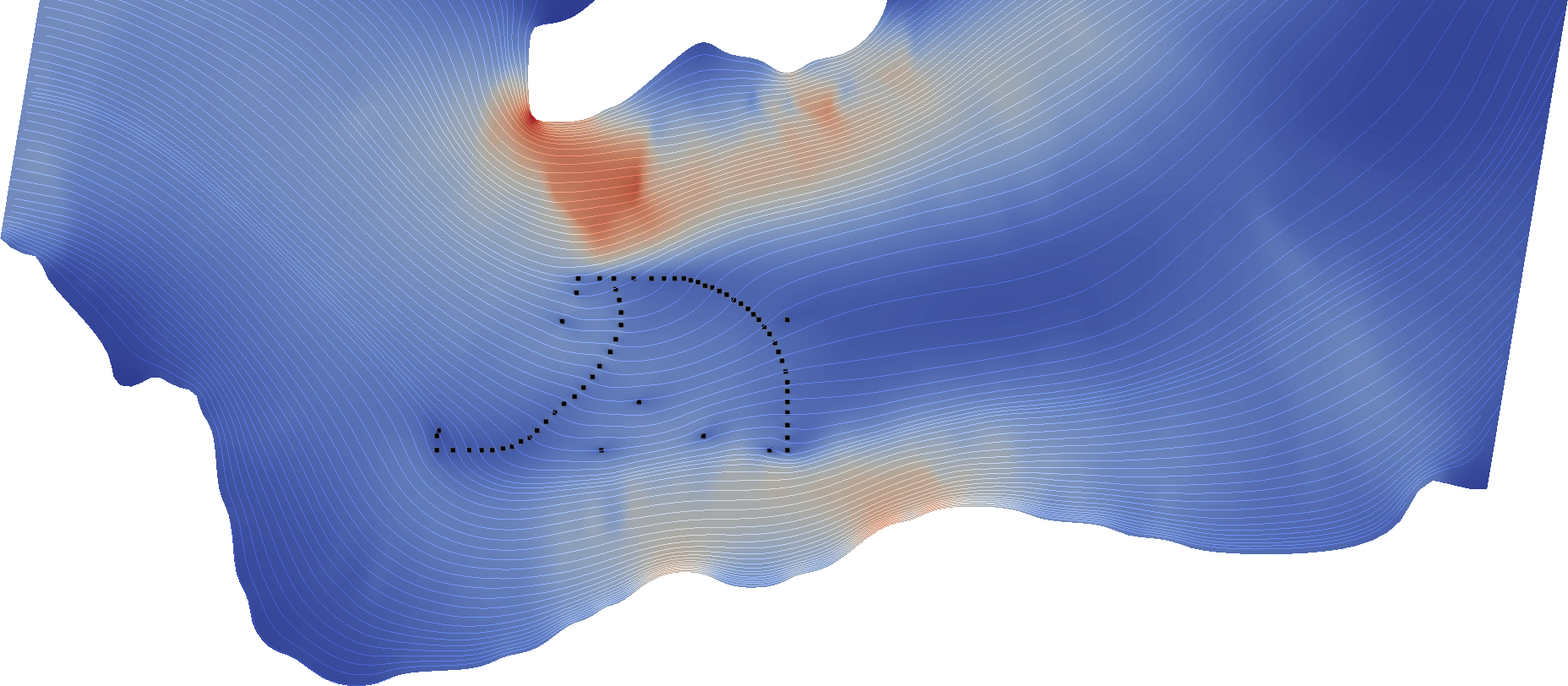}
  \caption{Streamline flow visualisation for southern part of the Pentland
    Firth scenario with the turbine layout acquired using the adjoint
    wake$\rightarrow$adjoint SWE optimisation for 64 turbines. Turbines
    toward the edge of each barrage lie further upstream than those in the
    centre which are approximately aligned orthogonal to the flow.
  }
\label{fig:orkney-streamlines}
\end{figure*}

Optimised layouts are displayed in \cref{fig:orkney-layouts}. For simplicity,
no minimum distance constraints (to enforce a minimum distance between
turbines) were used for optimisation. The layout produced by the genetic
algorithm features many overlapping turbines which the second stage adjoint
SWE optimisation is unable to rectify.  The layout which yields the greatest
power (`Adjoint~wake$\rightarrow$Adjoint~SWE') features two barrage-like
configurations of turbines. A streamline flow visualisation for this layout is
presented in \cref{fig:orkney-streamlines}. Both barrages have similar
features; towards the edges the turbines are positioned further upstream than
those close to the centre where they are approximately orthogonal to the flow
direction.  These features help to channel the flow through the array and are
conjectured to be the reason for enhanced performance of the design which
local optimisation alone was not able to achieve. \Citet{funke2014} optimise
for 128 and 256 turbines in the same domain but do not include the effects of
bathymetry or make use of global optimisation techniques. Whilst this study
uses fewer and larger turbines some comparisons may still be hypothesised.
Both cases feature barrages of turbines conjectured to guide the flow through
the array, however, the inclusion of bathymetry causes the curvature of the
barrages to be increased --- this is most likely to minimise the number of
turbines in the slowest part (south west corner) of the domain whilst still
guiding the flow toward the turbines in the faster flowing water in the north
eastern part of the turbine site. The large difference in flow velocities in
the turbine site is caused by the inclusion of bathymetry; \citet{funke2014}
use a constant depth of \SI{50}{\m} whilst in this study the bathymetry varies
between \SI{5}{\m} and \SI{105}{\m} with the depth in turbine site varying
from \SIrange{10.6}{35.5}{\m}. The large difference in depth between the north
and south of the domain and throughout the turbine site cause a much broader
range of flow magnitudes within the turbine site and is believed to be the
reason for the relocation of the barrages.

It is also noted that this layout does not share features with the first
optimisation stage (`Adjoint wake') which has approximately seven
south-south-west to north-north-east trending barrages of turbines (roughly
orthogonal to the flow direction). The difference in layouts is most likely
due to the simplicity in which the wakes are combined in the wake model.  An
improved wake combination method is more likely to result in a layout which
better matches that achieved after the second optimisation stage.  This would
also decrease the number of iterations required by the SWE model,
significantly reducing computation time. Whilst this relatively low resolution
scenario takes approximately \SI{6}{minutes} for each iteration of the SWE
model, the wake model takes approximately \SI{20}{seconds}. This saving would
naturally extend to higher resolution problems where computation time may be
significantly reduced.

\begin{table*}[tb]
  \makeatletter{}\newcommand{\raas}{$\rightarrow$Adjoint SWE}
\newcommand{\ra}{$\rightarrow$}
\begin{tabular}{c c c} \toprule
\textbf{Optimisation} & \textbf{Iterations}   & \textbf{Power (\si{\mega\W})} \\ \midrule
Initial               & --                    & \num{681.17}                  \\
Adjoint SWE           & \num{98}              & \num{730.49}                  \\
Adjoint wake\raas     & \num{864}, \num{130}  & \num{667.67}\ra\num{819.48}   \\
Genetic wake\raas     & \num{1000}, \num{103} & \num{643.96}\ra\num{729.64}   \\ \bottomrule
\end{tabular}

  \caption{Optimised power from the Pentland Firth scenario. The two stage
    adjoint optimisation of the wake model yields the highest power yet the
    power from the layout of the first stage is lower than the initial layout.
    The case is similar for the two stage genetic optimisation of the wake
    model but a second stage optimised power is similar to that of local
    optimisation of the SWE model.
  }
  \label{tab:orkney}
\end{table*}

\makeatletter{}\section{Conclusions}\label{sec:conclusion}

This work presents a wake model which has been demonstrated to act as a proxy
to the shallow water model in \otf \citep{funke2014}. The wake model is based
upon a reference wake pre-computed using a shallow water equation solver. It
provides a relatively cheap computation of the power extracted from a given
turbine array layout design and thus allows for global optimisation schemes
requiring large numbers of iterations to be used.

The speed of the wake model is its primary strength, however, it suffers in a
number of areas. Notably the model is not accurate when the turbine site
becomes saturated with turbines. Whether this is actually a problem in
realistic industrial applications is yet to be considered (an overly saturated
site may not be as viable as a lesser saturated site due to diminishing
marginal returns \citep{culley2014b}. The inaccuracy of the model is likely to
do with the simplistic manner in which wakes from multiple turbines are
combined and improvement to this should be considered in future work.
Improvements could also be made to the model by generating a reference wake
based on a flow speed which is more representative of the site in question
and scaling it appropriately.

Experiments utilising the wake model compared genetic and gradient-based
algorithms with basin-hopping (i.e.\ two global optimisers) and show that
layouts are generally optimised more effectively using gradient-based
approaches. This is particularly true when large numbers of turbines are
considered. However, optimisation via genetic algorithms show that seeding the
problem with a number of sensible guesses can significantly reduce the number
of iterations required for a solution to be found. Future work may extend this
idea by considering alternative methods of generating initial turbine layouts
as this may significantly decrease the number of iterations required during
optimisation.

When optimising layouts in a two-stage approach (global optimisation of the
wake model followed by local optimisation of the SWE model), the power
extracted from the array is similar to or greater than the power
achieved from just locally optimising the array using the SWE model, that is it
enables us to mitigate the problem of local maxima. In cases where bathymetry
is included, the improvements using two-stage optimisation are shown to be
significant with up to \SI{25}{\percent} improvements for idealised cases and
\SI{12}{\percent} for the more realistic Pentland Firth scenario. Whilst not
all cases see a reduction in computation time (due to the wake model producing
a poor layout for the second stage of optimisation) many see a reduction in
computation time of the order \SIrange{30}{40}{\percent}. Improvements in the
extracted power of this size could be the difference between a turbine array
being economically viable or not.

\makeatletter{}\section*{Acknowledgements}

The Authors would like to acknowledge financial support from Imperial College
London and the Engineering and Physical Sciences Research Council (grant
references: EP/L000407/1 and EP/J010065/1). Valuable discussions contributing
to this work were held with Stephan Kramer and Dave Culley.

\clearpage
\makeatletter{}\addcontentsline{toc}{section}{References}
\bibliography{sources}
\bibliographystyle{elsarticle-harv}

\end{document}